\input ./style/arxiv-general.cfg
\documentclass[aos,MSNbibl,nameyear,seceqn,dvips]{arximspdf}
\makeatletter
   \@ifpackageloaded{graphicx}{}{\usepackage{graphicx}}
\makeatother
\usepackage{mathrsfs}
\usepackage{mathbh}

\doi{10.1214/15-AOS1366}
\volume{44}
\issue{1}
\pubyear{2016}
\firstpage{255}
\lastpage{287}
\docsubty{FLA}

\makeatletter
\newcommand{\rrvert}{\vert}
\newcommand{\llvert}{\vert}
\newcommand{\eqref}[1]{(\ref{#1})}
\newtheorem{thmm}{Theorem}[section]
\newtheorem{lemma}[thmm]{Lemma}
\newproclaim{definition}[thmm]{Definition}

\newproclaim{remark}[thmm]{Remark}

\newcommand{\Var}{\operatorname{Var}}
\newcommand{\E}{\mathrm{E}}

\def\E{\mathbb{E}}
\def\N{\mathbb{N}}
\def\R{\mathbb{R}}
\def\N{\mathbb{N}}
\def\P{\mathbb{P}}
\def\bblambda{\bolds{\lambda}}

\newcommand{\barbeta}{\stackrel{}{\bar{\beta}}}
\newcommand{\Leb}{d \bblambda^d}
\makeatother

\begin{document}
\begin{frontmatter}

\title{Adaptation to lowest density regions with application to support
recovery\thanksref{T1}}
\runtitle{Adaptation to lowest density regions}
\thankstext{T1}{Supported by the DFG Priority Program SPP 1324, RO 3766/2-1.}

\begin{aug}
\author[A]{\fnms{Tim}~\snm{Patschkowski}\corref{}\ead[label=e1]{tim.patschkowski@ruhr-uni-bochum.de}}
\and
\author[A]{\fnms{Angelika} \snm{Rohde}\ead[label=e2]{angelika.rohde@ruhr-uni-bochum.de}}
\runauthor{T. Patschkowski and A. Rohde}
\affiliation{Ruhr-Universit\"at Bochum}

\address[A]{Fakult\"at f\"ur Mathematik\\
Ruhr-Universit\"at Bochum\\
44780 Bochum\\
Germany\\
\printead{e1}\\
\phantom{E-mail:\ }\printead*{e2}}

\end{aug}

%
\received{\smonth{9} \syear{2014}}
%
\revised{\smonth{5} \syear{2015}}

%
\begin{abstract}
A scheme for locally adaptive bandwidth selection is proposed which
sensitively shrinks the bandwidth of a kernel estimator at lowest
density regions such as the support boundary which are unknown to the
statistician.
In case of a H\"older continuous density, this locally minimax-optimal
bandwidth is shown to be smaller than the usual rate, even in case of
homogeneous smoothness. Some new type of risk bound with respect to a
density-dependent standardized loss of this estimator is established.
This bound is fully nonasymptotic and allows to deduce convergence
rates at lowest density regions that can be substantially faster than
$n^{-1/2}$. It is complemented by a weighted minimax lower bound which
splits into two regimes depending on the value of the density. The new
estimator adapts into the second regime, and it is shown that
simultaneous adaptation into the fastest regime is not possible in
principle as long as the H\"older exponent is unknown. Consequences on
plug-in rules for support recovery are worked out in detail. In
contrast to those with classical density estimators, the plug-in rules
based on the new construction are minimax-optimal, up to some
logarithmic factor.
\end{abstract}

%
\end{frontmatter}

\section{Introduction}\label{sec1}
Adaptation in the classical context of nonparametric function
estimation in Gaussian white noise has been extensively studied in the
statistical literature. Since \citet{nussbaum1996} has established
asymptotic equivalence in Le Cam's sense for the nonparametric models
of density estimation and Gaussian white noise, a rigorous framework is
provided which allows to carry over specific statistical results
established for the Gaussian white noise model to the model of density
estimation, at least in dimension one. Density estimation is as one of
the most fundamental problems in statistics subject to a variety of
recent studies; see, for example, \citet{efromovich2008},
\citet{gachnicklspokoiny2013},
\citet{lepski2013}, \citet{birge2014} and \citet
{liuwong2014}. It has become clear that under the conditions for the
asymptotic equivalence to hold, minimax rates of convergence in density
estimation with respect to pointwise or mean integrated squared error
loss coincide with the optimal convergence rates obtained in the
context of nonparametric regression, and the procedures are typically
identical on the level of ideas. A main requisite on the density for
Nussbaum's (\citeyear{nussbaum1996}) asymptotic equivalence is the
assumption that it is
compactly supported and uniformly bounded away from zero on its
support. If this assumption is violated, the density estimation
experiment may produce statistical features which do not have any
analog in the regression context. For instance, minimax estimation of
noncompactly supported densities under $L_p$-loss bears striking
differences to the compact case; see \citet
{judiskylambert-lacroix2004}, \citet
{reynaud-bouretrivoirardtuleau-malot2011} and
Goldenshluger and Lepski (\citeyear
{goldenshlugerlepski2011,goldenshlugerlepski2013}). The minimax rates
reflect an
interplay of the regularity parameters and the parameter of the loss
function, an effect which is caused by the tail behavior of the
densities under consideration. In this article, we recover such an
exclusive effect even for compactly supported densities. It turns out
that minimax estimation in regions where the density is small is
possible with higher accuracy although fewer observations are
available, leading to rates which can be substantially faster than
$n^{-1/2}$. Even more, this accuracy can be achieved to a large extent
without a priori knowledge of these regions by a kernel density
estimator with an adaptively selected bandwidth. As discovered by
\citet{butucea2001}, the exact constant of normalization for pointwise
adaptive univariate density estimation on Sobolev classes depends
increasingly on the density at the point of estimation itself. The
crucial observation is that the classical bias variance trade-off does
not reflect the dependence of the kernel estimator's variance on the
density, which brings the idea of an estimated variance in the
bandwidth selection rule into play. Although Butucea's interesting
result requires the point of estimation to be fixed, it suggests that a
potential gain in the rate might be possible at lowest density regions.
In this paper, we investigate the problem of adaptation to lowest
density regions under anisotropic H\"older constraints. A bandwidth
selection rule is introduced which provably attains fast pointwise
rates of convergence at lowest density regions. On this way, new
weighted lower risk bounds over anisotropic H\"older classes are
established, which split into two regimes depending on the value of the
density. We show that the new estimator uniformly improves the global
minimax rate of convergence, adapts to the second regime and finally
that adaptation into the fastest regime is not possible in principle if
the density's regularity is unknown. We identify the best possible
adaptive rate of convergence
\[
n^{-{\barbeta}/{(\barbeta+d)}}
\]
(up to a logarithmic factor), where $\barbeta$ is the unnormalized
harmonic mean of the $d$-dimensional H\"older exponent.

This breakpoint determines the attainable speed of convergence of
plug-in estimators for functionals of the density where the quality of
estimation at the boundary is crucial. We exemplarily demonstrate it
for the problem of support recovery. In order to line up with the
related results of \citet{cuevasfraiman1997} about plug-in rules for
support estimation and \citet{rigolletvert2009} on minimax
analysis of
plug-in level-set estimators, we measure the performance of the plug-in
support estimator with respect to the global measure of symmetric
difference of sets under the margin condition [\citet
{polonik1995}; see also
\citet{mammentsybakov1999} and \citet{tsybakov2004}].
In contrast to level set estimation, however, plug-in rules for the
support functional possess sub-optimal convergence rates when the
classical kernel density estimator with minimax-optimal global
bandwidth choice is used. We determine the optimal minimax rate for
support recovery
\[
n^{-{\gamma\beta}/{(\beta+d)}}
\]
(up to a logarithmic factor), where $\gamma$ denotes the margin
exponent, $d$ the dimension and $\beta$ the isotropic H\"older
exponent. Our result demonstrates that support recovery is possible
with higher accuracy than level set estimation as already conjectured
by \citet{tsybakov1997}. We finally show that the performance of
the plug-in
support estimator resulting from our new density estimator turns out to
be minimax-optimal up to a logarithmic factor.

The article is organized as follows. Section~\ref{sec:preliminairies}
contains the basic notation. In Section~\ref{sec:densitybounds}, the
adaptive density estimator is introduced, new
weighted lower pointwise risk bounds are derived and the optimality
performance of the estimator is proved. Section~\ref
{sec:supportapplication} addresses the important problem of density
support estimation as an example of a functional which substantially
benefits from the new density estimator.
The proofs are deferred to Section~\ref{sec:proofs} and the
supplemental article [\citet{supp}].

\section{Preliminaries and notation} \label{sec:preliminairies}

All our estimation procedures are based on a sample of $n$ real-valued
$d$-dimensional random vectors $X_i = (X_{i,1}, \ldots,\break X_{i,d})$,
$i = 1, \ldots, n$ ($d \geq1$ and if not stated otherwise $n \geq
2$), that are independent and identically distributed according to some
unknown probability measure $\mathbb{P}$ on $\R^d$ with continuous
Lebesgue density $p$. $\E_p^{\otimes n}$ denotes the expectation with
respect to the $n$-fold product measure $\mathbb{P}^{\otimes n}$. Let
\[
\hat{p}_{n,h}(t) = \hat{p}_{n,h}(t, X_1,
\ldots, X_n):= \frac{1}{n} \sum_{i=1}^n
K_h(t - X_i),
\]
denote the kernel density estimator with $d$-dimensional bandwidth $h =
(h_1,\break \ldots,   h_d)$ at point $t\in\R^d$, where
\[
K_h(x) := \Biggl( \prod_{i=1}^d
h_i \Biggr)^{-1} K \biggl( \frac
{x_1}{h_1}, \ldots,
\frac{x_d}{h_d} \biggr)
\]
describes a rescaled kernel supported on $\prod_{i=1}^d [-h_i,h_i]$.
The kernel function $K$ is assumed to be compactly supported on
$[-1,1]^d$ and to be of product structure, that is, $K(x_1,\ldots,x_d)
= \prod_{i=1}^d K_i(x_i)$. Additionally, $K_{i,{h_i}}(x):= h_i^{-1}
K_i(x/h_i), i=1,\ldots,d$. The components $K_i$ are assumed to
integrate to one and to be continuous on its support with $K_i(0)>0$.
If not stated otherwise, they are symmetric and nonnegative, implying
that the kernel is of first order. Recall that $K$ is said to be of
$k$th order, $k=(k_1,\ldots,k_d)\in\N^d$, if the functions $x \mapsto
x_i^{j_i} K_i(x_i)$, $j_i\in\N$ with $1\leq j_i\leq k_i$, $i=1,\ldots
,d$, satisfy
\[
\int x_i^{j_i} K_i(x_i) \,d
\bblambda(x_i) = 0,
\]
where $\bblambda^d$ denotes the Lebesgue measure on $\mathbb{R}^d$
throughout the article. The Lebesgue measure on $\mathbb{R}$ is denoted
by $\bblambda$.
For any function $f: \mathbb{R}^d \rightarrow\mathbb{R}$ and $x =
(x_1, \ldots, x_d) \in\mathbb{R}^d$, we define the univariate
functions
%
%
\begin{eqnarray}
\label{fastfertig2} %
f_{i,x}: \mathbb{R} &\longrightarrow&\mathbb{R}
\nonumber
\\[-8pt]
\\[-8pt]
\nonumber
y & \longmapsto& f(x_1, \ldots, x_{i-1}, y,
x_{i+1}, \ldots, x_d)
\end{eqnarray}
and denote by $P_{y,l}^{(f_{i,x})}$ the Taylor polynomial
%
%
\begin{equation}
\label{taylorpolynomial} P_{y,l}^{(f_{i,x})}(\cdot) := \sum
_{k=0}^l \frac{f_{i,x}^{(k)}(y)}{k!} (\cdot-y)^k
\end{equation}
of $f_{i,x}$ at the point $y \in\mathbb{R}$ of degree $l$ (whenever it
exists). Let $\mathscr{H}_d(\beta,L)$ be the anisotropic H\"older class
with regularity parameters $(\beta,L)$, that is, any function $f$
belonging to this class fulfills for all $y,y' \in\mathbb{R}$ the inequality
\[
\sup_{x \in\mathbb{R}^d} \bigl\vert f_{i,x} (y) - f_{i,x}
\bigl(y^\prime\bigr) \bigr\vert\leq L \bigl\vert y - y^\prime
\bigr\vert^{\beta_i}
\]
for those $i \in\{ 1, \ldots, d \}$ with $\beta_i \leq1$, and in case
$\beta_i > 1$ admits derivates with respect to its $i$th coordinate up
to the order $\lfloor\beta_i \rfloor:= \max\{ n \in\mathbb{N} : n <
\beta_i \}$, such that the approximation by the Taylor polynomial satisfies
\[
\sup_{x \in\mathbb{R}^d} \bigl\llvert f_{i,x}(y) - P_{y',\lfloor\beta
_i \rfloor}^{(f_{i,x})}(y)
\bigr\rrvert\leq L \bigl\vert y - y'\bigr \vert^{\beta_i}\qquad \mbox{for
all } y, y' \in\mathbb{R}.
\]
For adaptation issues, it is assumed that $\beta= (\beta_1, \ldots,
\beta_d) \in\prod_{i=1}^d [ {\beta}^*_{i,l}, {\beta}_{i,u}^* ]$
and $L
\in[{L}_l^*, {L}_u^*]$ for some positive constants ${\beta}^*_{i,l} <
{\beta}_{i,u}^*$, $i=1,\ldots,d$, and ${L}_l^* < {L}_u^*$. For short, we
simply write $\beta^*$ and $L^*$ for the couples $(\beta^*_l,\beta
^*_u)$ and $(L^*_l,L^*_u)$, and finally $\mathcal{R}(\beta^*,L^*)$ for
the rectangle $\prod_{i=1}^d[ {\beta}^*_{i,l}, {\beta
}_{i,u}^*]\times
[{L}_l^*, {L}_u^*]$. It turns out that all rates of convergence
emerging in an anisotropic setting involve the unnormalized harmonic
mean of the smoothness parameters
\[
\barbeta:= \Biggl( \sum_{i=1}^d
\frac{1}{\beta_i} \Biggr)^{-1}.
\]
To focus on rates only and for ease of notation, we denote by $c$
positive constants that may change from line to line. All relevant
constants will be numbered consecutively. Dependencies of the constants
on the functional classes' parameters are always indicated and it
should be kept in mind that the constants can potentially depend on the
chosen kernel, the loss function and the dimension as well.
Furthermore, $\mathscr{P}_d(\beta,L)$ denotes the set of all
probability densities in $\mathscr{H}_d(\beta,L)$. It is well known
that any function $f \in\mathscr{P}_d(\beta,L)$ is uniformly bounded
by a constant
%
%
\begin{equation}
\label{c1} c_1(\beta,L) = \sup\bigl\{ \Vert p
\Vert_{\mathrm{sup}} : p \in\mathscr{P}_d(\beta,L) \bigr\}
\end{equation}
depending on the regularity parameters only.

\section{New lower risk bounds, adaptation to lowest density regions}
\label{sec:densitybounds}
The fully nonparametric problem of estimating a density $p$ at some
given point $t = (t_1, \ldots, t_d)$ has quite a long history in the
statistical literature and has been extensively studied. Considering
different estimators, a very natural question is whether there is an
estimator that is optimal and how optimality can be exactly described.
A common concept of optimality is stated in a minimax framework. An
estimator $T_n(t) = T_n(t,X_1, \ldots, X_n)$ is called minimax-optimal
over the class $\mathscr{P}_d(\beta,L)$ if its risk matches the
minimax risk
\[
\inf_{T_n(t)} \sup_{p \in\mathscr{P}_d(\beta
,L)} \E_p^{\otimes n}
\bigl\vert T_n(t) - p(t)\bigr \vert^r
\]
for some $r \geq1$, where the infimum is taken over all estimators.
However, the minimax approach is often rated as quite pessimistic as it
aims at finding an estimator which performs best in the worst
situation. Different in spirit is the oracle approach. Within a
pre-specified class $\mathscr{T}$ of estimators, it aims at finding for
any individual density the estimator $\hat{T}_n\in\mathscr{T}$ which is
optimal, leading to oracle inequalities of the form
\[
\E_p^{\otimes n} \bigl\vert\hat{T}_n(t)- p(t)
\bigr\vert^r \leq c \inf_{T_n\in\mathscr{T}} \E_p^{\otimes n}
\bigl\vert T_n(t) - p(t) \bigr\vert^r + R_n(t)
\]
with a remainder term $R_n(t)$ depending on the class $\mathscr{T}$,
the underlying density $p$ and the sample size only. Besides having the
drawback that there is no notion of optimality judging about the
adequateness of the estimator's class, an equally severe problem may be
caused by the fact that the remainder term is uniform in $\mathscr{T,}$
and thus a worst case remainder. The latter is responsible for the fact
that our fast convergence rates cannot be deduced from the oracle
inequality in Goldenshluger and \citet{lepski2013}, the order for their
remainder being unimprovable, however. In this article, we introduce
the notion of best possible $p$-dependent minimax speed of convergence
$\psi^n_{p(t),\beta,L}$ within the function class\vspace*{1pt} $\mathscr
{P}_d(\beta,
L)$ and aim at constructing an estimator $T_n(t)$ bounding the risk
\[
\sup_{p \in\mathscr{P}_d(\beta,L)} \mathop{\sup_{t
\in\mathbb{R}^d:}}_{ p(t)>0}
\E_p^{\otimes n} \biggl( \frac{
\vert T_n(t) - p(t) \vert}{ \psi^n_{p(t),\beta,L} } \biggr)^r
\]
uniformly over a range of parameters $(\beta,L)$. First, this requires
a suitable definition of the quantity $\psi^n_{p(t),\beta,L}$.

\subsection{New weighted lower risk bound} As we want to work out the
explicit dependence on the value of the density, it seems suitable to
fix an arbitrary constant $\varepsilon\in(0,1)$, and to pick out
maximal not necessarily disjoint subsets $U_{\delta}$ of $\mathscr
{P}_d(\beta,L)$ with the following properties: $\cup U_{\delta}=\{
p\in
\mathscr{P}_d(\beta,L): p(t)>0\}$, and pairwise ratios $p(t)/q(t)$,
$p,q\in U_{\delta}$, are bounded away from zero by $\varepsilon$ and
from infinity by $1/\varepsilon$. This motivates the construction of
the subsequent theorem.

%
\begin{thmm}[(New weighted lower risk bound)] \label{lowerbounddensity}
For any $\beta=(\beta_1,\ldots,\beta_d)$ with $0 < \beta_i \leq
2$, $i
= 1, \ldots, d$, $L>0$ and $r \geq1$, there exist constants\break
$c_2(\beta
,L,r)>0$ and $n_0(\beta,L) \in\mathbb{N}$, such that for every $t\in
\R
$ the pointwise minimax risk over H\"older-smooth densities is bounded
from below by
\[
\inf_{0 < \delta\leq c_1(\beta,L)} \inf_{T_n(t)} \mathop{\sup_{p \in\mathscr{P}_d(\beta, L): }}_{ \delta/2
\leq
p(t) \leq\delta} \E_p^{\otimes n} \biggl(
\frac{\vert
T_n(t)-p(t) \vert}{\psi_{p(t),\beta}^n} \biggr)^r \geq c_2(\beta,L,r)
\]
for all $n \geq n_0(\beta,L)$, where $\psi_{x,\beta}^n:=x\wedge
(x/n)^{{\barbeta}/{(2\barbeta+1)}}$ and $c_1(\beta,L)$ defined in~\eqref{c1}.
\end{thmm}

%
\begin{remark}
(i) The lower bound of the above theorem is attained by the oracle
estimator
%
%
\begin{equation}
\label{oracle} T_n(t) := \hat{p}_{n,h_{n,\delta}}(t) \cdot\mathbh{1}
\bigl\{ \delta\geq n^{-\barbeta/(\barbeta+1)} \bigr\}
\end{equation}
with $h_{n,\delta,i}=(\delta/n)^{({1}/{(2\barbeta+1)}) ({1}/{\beta_i})}$. Hence, $\psi_{p(t),\beta}^n$ cannot be improved in principle. We
refer to it in the sequel as $p$-dependent speed of convergence within
the functional class $\mathscr{P}_d(\beta,L)$.

%
\begin{figure}[b]

\includegraphics{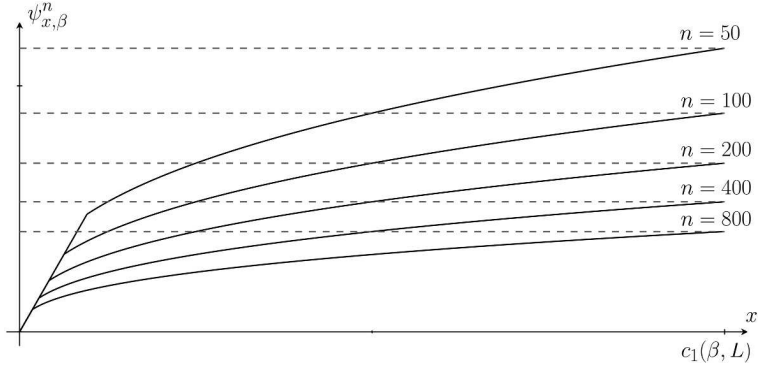}

\caption{New lower bound (solid line), classical lower bound (dashed line).}
\label{fig:lowerboundsplit}
\end{figure}

(ii) Note that for the classical minimax rate $n^{-\barbeta
/(2\barbeta+1)}$,
\[
\lim_{n \rightarrow\infty} \inf_{0 < \delta\leq
c_1(\beta,L)} \inf_{T_n(t)}
\mathop{\sup_{p
\in\mathscr{P}_d(\beta, L):}}_{ \delta/2 \leq p(t) \leq\delta
} \E_p^{\otimes n}
\biggl(\frac{\vert T_n(t)-p(t) \vert
}{n^{-\barbeta/(2\barbeta+1)}} \biggr)^r = 0
\]
as a direct consequence of the subsequently formulated Theorem~\ref
{upperbounddensity}. The $p$-dependent speed of convergence $\psi
_{p(t),\beta}^n$ is of substantially smaller order than the classical one
along a shrinking neighborhood of lowest density regions.
\end{remark}

Note that the exponent $\barbeta/(2 \barbeta+ 1)$ implicitly depends
on the dimension $d$ and coincides in case of isotropic smoothness with
the well-known exponent $\beta/(2 \beta+ d)$. It splits into two
regimes which are listed and specified in the following table.\vspace*{6pt}

%
%
\noindent
\begin{tabular*}{\textwidth}{@{\extracolsep{\fill}}lcc@{}}
\hline
& {\fontsize{9}{11}\selectfont{\textbf{Regime}}} &
\multicolumn{1}{c@{}}{{\fontsize{9}{11}\selectfont{\textbf{Rate}
$\bolds{\psi_{x,\beta}^n}$}}}\\
\hline
\phantom{i}(i) & $x \leq n^{-{\barbeta}/{(\barbeta+1)}}$ & $x$ \\
(ii) & $n^{-{\barbeta}/{(\barbeta+1)}} < x \leq c_1(\beta,L)$ &
$ ( \frac{x}{n} )^{{\barbeta}/{(2\barbeta+1)}}$\\
\hline
\end{tabular*}\vspace*{9pt}\label{table1}

The worst $p$-dependent speed of convergence within $\mathscr
{P}_d(\beta
,L)$, namely
\[
\sup_{0<x\leq c_1(\beta,L)}\psi_{x,\beta}^n,
\]
reveals the classical minimax rate $n^{-\barbeta/(2\barbeta+1)}$. The
fastest rate in regime (ii) is of the order
\[
n^{-\barbeta/(\barbeta+1)} \qquad\mbox{for } x=n^{-\barbeta
/(\barbeta+1)},
\]
which is substantially smaller than the classical minimax risk bound.
Figure~\ref{fig:lowerboundsplit} visualizes the split-up into the
regimes and relates the new $p$-dependent rate of Theorem~\ref
{lowerbounddensity} to the classical minimax rate for different sample
sizes from $n = 50$ to $n = 800$.

It becomes apparent from the proof that the lower bound actually even
holds for the subset of $(\beta,L)$-regular densities with compact
support. At first glance, however, the new lower bound is of
theoretical value only, because the value of a density at some point to
be estimated is unknown. The question is whether it is possible to
improve the local rate of convergence of an estimator without prior
knowledge in regions where fewer observations are available, that is,
to which extent it is possible to adapt to lowest density regions.

\subsection{Adaptation to lowest density regions} \label{construction!}
Adaptation is an important challenge in nonparametric estimation.
Lepski (\citeyear{lepski1990}) introduced a sequential multiple testing
procedure for
bandwidth selection of kernel estimators in the Gaussian white noise
model. It has been widely used and refined for a variety of adaptation
issues over the last two decades. For recent references, see \citet
{ginenickl2010}, \citet{chichignoud2012},
Goldenshluger and Lepski (\citeyear
{goldenshlugerlepski2011,goldenshlugerlepski2013}), \citet
{chichignoudlederer2014}, Jirak, Meister and Rei{\ss}
(\citeyear{jirakmeisterreiss2014}), \citet
{dattnerreisstrabs2014} and \citet{bertinlacourrivoirard2014}
and \citet{lepski2014}
among many others.
Our subsequently constructed estimator is based on the anisotropic
bandwidth selection procedure of \citet
{kerkyacharianlepskipicard2001}, which has been developed in the
Gaussian white noise model,
but incorporates the new approach of adaptation to lowest density
regions. Although Goldenshluger and \citet{lepski2013} pursue a
similar goal
via some kind of empirical risk minimization, their oracle inequality
provides no faster rates than $n^{-1/2}$ times the average of the
density over the unit cube around the point under consideration. They
deduce from it adaptive minimax rates of convergence with respect to
the $L_p$-risk over anisotropic Nikol'skii classes for density
estimation on $\R^d$. As concerns adaptation to lowest density regions
such as the unknown support boundary, this oracle inequality is not
sufficient as no faster rates than $n^{-1/2}$ can be deduced from it,
and it is not clear whether these faster rates are attainable for their
estimator in principle. Besides having the drawback that there is no
notion of optimality judging about the adequateness of the estimator's
class, an equally severe problem of the oracle approach may be caused
by the fact that the remainder term is uniform in the estimator's
class, and thus a worst case remainder. The latter is responsible for
the fact that our fast convergence rates cannot be deduced from the
oracle inequality in Goldenshluger and \citet{lepski2013}, the
order for
their remainder being unimprovable, however. It raises the question
whether this imposes a fundamental limit on the possible range of
adaptation (the corresponding inequality resulting from the bound on
$\P^{\otimes n}(B_{1,m})$ has to be satisfied as well).
We shall demonstrate in what follows that it is even
possible to attain substantially faster rates, indeed that adaptation
to the whole second regime of Theorem~\ref{lowerbounddensity} is an
achievable goal, and that this describes precisely the full range where
adaptation to lowest density regions is possible as long as the
density's regularity is unknown.
Our procedure uses kernel density estimators $\hat{p}_{n,h}(t)$ with
multivariate bandwidths $h = (h_1, \ldots, h_d)$, which are able to
deal with different degrees of smoothness in different coordinate
directions. Note that optimal bandwidths for estimation of H\"
older-continuous densities are typically derived by a bias-variance
trade-off balancing the bias bound
%
%
\begin{equation}
\label{bias!} \bigl\llvert p(t) - \E_p^{\otimes n}
\hat{p}_{n,h}(t) \bigr\rrvert\leq c(\beta,L) \cdot\sum
_{i=1}^d h_i^{\beta_i},
\end{equation}
see \eqref{bias} in Section~\ref{sec:proofs} for details, against the
rough variance bound
%
%
\begin{equation}
\Var\bigl(\hat{p}_{n,h}(t) \bigr) \leq\frac{c_1(\beta,L) \Vert K \Vert
_2^2}{n \prod_{i=1}^{d} h_i},
\label{roughvariancebound}
\end{equation}
where $\Vert\cdot\Vert_2$ is the Euclidean norm [on $L_2(\bblambda
^d)$]. This bound leads to suboptimal rates of convergence whenever the
density is small since it is not able to capture small values of $p$ in
a small neighborhood around $t$
in contrast to the sharp convolution bound
%
%
\begin{equation}
\label{varianz!} \Var\bigl(\hat{p}_{n,h}(t) \bigr) \leq
\frac{1}{n} \bigl((K_h)^2 \ast p \bigr) (t) =:
\sigma_t^2(h).
\end{equation}
Balancing \eqref{bias!} and \eqref{varianz!} leads to smaller
bandwidths at lowest density regions as compared to bandwidths
resulting from the classical bias-variance trade-off between \eqref
{bias!} and \eqref{roughvariancebound}.
The convolution bound \eqref{varianz!} is unknown and it is natural to
replace it by its unbiased empirical version
\[
\tilde{\sigma}_t^2(h) := \frac{1}{n^2 \prod_{i=1}^d h_i^2} \sum
_{i=1}^n K^2 \biggl(
\frac{t-X_i}{h} \biggr).
\]
However, $\tilde{\sigma}_t^2(h)$ concentrates extremely poorly around
its mean if the bandwidth $h$ is small, which is just the important
situation at lowest density regions. Precisely, Bernstein's inequality
provides the bound
%
%
\begin{equation}\qquad
\label{bern} \mathbb{P}^{\otimes n} \biggl( \biggl\vert\frac{\tilde
{\sigma
}_{t}^2(h)}{\sigma_{t}^2(h)}-1
\biggr\vert\geq\eta\biggr) \leq2 \exp\Biggl( -\frac{3 \eta
^2}{2(3+2\eta) \Vert K \Vert_\mathrm{sup}^2}
\sigma_t^2(h) \cdot n^2\prod
_{i=1}^dh_i^2 \Biggr),
\end{equation}
which suggests to study the following truncated versions instead:
%
%
\begin{eqnarray}
\label{sigma_trunc} %
\sigma_{t,\mathrm{trunc}}^2(h) &:= &\max
\biggl\{ \frac{\log^2 n}{n^2
\prod_{i=1}^d h_i^2}, \sigma_t^2(h) \biggr\},
\nonumber
\\[-8pt]
\\[-8pt]
\nonumber
\tilde{\sigma}_{t,\mathrm{trunc}}^2(h) &:=& \max\biggl\{
\frac{\log^2
n}{n^2 \prod_{i=1}^d h_i^2}, \tilde{\sigma}_t^2(h) \biggr\}.
\end{eqnarray}
Without the logarithmic term, the truncation level ensures tightness of
the family of random variables $\tilde{\sigma}_{t,\mathrm
{trunc}}^2(h)/\sigma_{t,\mathrm{trunc}}^2(h)$, because the exponent in
\eqref{bern} remains a nondegenerate function in $\eta$. The
logarithmic term is introduced in order to guarantee sufficient
concentration of $\sup_{h}\vert1- \tilde{\sigma}_{t,\mathrm
{trunc}}^2(h)/\sigma_{t,\mathrm{trunc}}^2(h)\vert$.

\textit{Construction of the adaptive estimator.} Our estimation
procedure is developed in the anisotropic setting, in which neither
the variance bound nor the bias bound provides an immediate monotone
behavior in the bandwidth. Unlike in the univariate or isotropic
multivariate case, Lepski's (\citeyear{lepski1990}) idea of mimicking
the bias-variance
trade-off fails. Consequently, our estimation scheme imitates the
anisotropic procedure of \citet{kerkyacharianlepskipicard2001}
and \citet{klutchnikoff2005}, developed in the Gaussian white
noise model,
with the following changes. First, their threshold given by the
variance bound in the Gaussian white noise setting is replaced
essentially with the truncated estimate in \eqref{sigma_trunc}, which
is sensitive to small values of the density. Moreover, it is crucial in
the anisotropic setting that our procedure uses an ordering of
bandwidths according to these estimated variances instead of an
ordering according to the product of the bandwidth's components. The
bandwidth selection scheme chooses a bandwidth in the set
\[
\mathcal{H} := \Biggl\{ h = (h_1, \ldots,
h_d) \in\prod_{i=1}^d
(0,h_{\mathrm{max},i}] : \prod_{i=1}^d
h_i \geq\frac{\log^2 n}{n} \Biggr\},
\]
where for simplicity we set $(h_{\mathrm{max},1}, \ldots, h_{\mathrm
{max},d}) = (1, \ldots, 1)$. Let furthermore
\[
\mathcal{J} := \Biggl\{ j = (j_1, \ldots, j_d) \in
\mathbb{N}_0^d : \sum_{i=1}^d
j_i \leq\biggl\lfloor\log_2 \biggl(\frac{n}{\log^2
n}
\biggr) \biggr\rfloor\Biggr\}
\]
be a set of indices and denote by
\[
\mathcal{G} := \bigl\{ \bigl(2^{-j_1}, \ldots, 2^{-j_d} \bigr) :
j \in\mathcal{J} \bigr\} \subset\mathcal{H}
\]
the corresponding dyadic grid of bandwidths, that serves as a
discretization for the multiple testing problem in Lepski's selection
rule. For ease of notation, we abbreviate dependences on the bandwidth
$(2^{-j_1}, \ldots, 2^{-j_d})$ by the multi-index $j$. Next,
with $j \wedge m$ denoting the minimum by component,
the set of admissible bandwidths is defined as
%
%
\begin{eqnarray}
\label{admissiblebandwidths} %
\mathcal{A} &= &\mathcal{A}(t)\nonumber\\
 &:=& \Bigl\{ j \in
\mathcal{J} :\bigl \vert\hat{p}_{n,j \wedge m} (t) - \hat{p}_{n,m}(t)
\bigr\vert\leq c_3 \sqrt{\hat{\sigma}_t^2(m)
\log n}
\\
&&{}\mbox{for all } m \in\mathcal{J} \mbox{ with } \hat{\sigma}_t^2(m)
\geq\hat{\sigma}_t^2(j) \Bigr\},\nonumber
\end{eqnarray}
with a properly chosen constant
$c_3 = c_3(\beta^*,L^*)$ satisfying the constraint \eqref{c4_2}
appearing in the proof of Theorem~\ref{upperbounddensity}.
Here, both the threshold and the ordering of bandwidths are defined via
the truncated variance estimator
%
%
\begin{eqnarray}
\label{truncated_variance_estimator} 
\hat{\sigma}_t^2(h) &:= &\min\biggl
\{
\tilde{\sigma}_{t,\mathrm
{trunc}}^2(h), \frac{\Vert K \Vert_2^2 c_1}{n \prod_{i=1}^d h_i} \biggr
\}
\nonumber
\\[-8pt]
\\[-8pt]
\nonumber
&= &\min\Biggl\{ \max\Biggl[ \frac{\log^2 n}{n^2 \prod_{i=1}^d h_i^2},
\frac{1}{n^2 \prod_{i=1}^d h_i^2} \sum
_{i=1}^n K^2 \biggl(
\frac
{t-X_i}{h} \biggr) \Biggr], \frac{\Vert K \Vert_2^2 c_1}{n \prod
_{i=1}^d h_i} \Biggr\},\hspace*{-25pt}
\end{eqnarray}
where $c_1 = c_1(\beta^*, L^*)$ is an upper bound on $c_1(\beta,L)$ in
the range of adaptation.
The threshold in \eqref{admissiblebandwidths} could be modified by a
further logarithmic factor to avoid the dependence of the constants on
the range of adaptation.
Recall again that this refined estimated threshold is crucial for our
estimation scheme.
The procedure selects the bandwidth among all admissible bandwidths with
%
%
\begin{equation}
\label{lepski_minimization} \hat{j} = \hat{j}(t) \in\mathop{\arg\min}_{j \in
\mathcal{A}} \hat{
\sigma}_t^2(j).
\end{equation}
Finally,
\[
\hat{p}_n := \hat{p}_{n,\hat{j}} \wedge c_1
\]
defines the adaptive estimator. In case of isotropic H\"older
smoothness, it is sufficient to restrict the grid to bandwidths with
equal components, and we even simplify the method by replacing the
ordering by estimated variances in condition \eqref
{truncated_variance_estimator}
''for all $ m \in\mathcal{J} \mbox{ with } \hat{\sigma}_t^2(m)
\geq
\hat{\sigma}_t^2(j)$'' by the classical order ``for all $m \in
\mathcal
{J} \mbox{ with } m\geq j$'' as the componentwise ordering is the same
for all components.

\textit{Performance of the adaptive estimator.} Clearly, the
truncation in the threshold imposes serious limitations to which extent
adaptation to lowest densities regions is possible. However, a careful
analysis of the ratio
\[
\sup_h \biggl\vert\frac{\tilde{\sigma}_{t,\mathrm
{trunc}}^2(h)}{\sigma_{t, \mathrm{trunc}}^2(h)}-1 \biggr\vert
\]
rather than the difference $\sup_h\vert\tilde{\sigma
}_{t,\mathrm
{trunc}}^2(h) - \sigma_{t, \mathrm{trunc}}^2(h)\vert$ allows to
prove indeed that adaptation is possible in the whole second regime.

%
\begin{thmm}[(New upper bound)] \label{upperbounddensity}
For any rectangle $\mathcal{R}(\beta^*,L^*)$ with $[\beta
^*_{i,l},\beta
^{*}_{i,u}]\subset(0,2]$, $[L^*_{l},L^*_{u}]\subset(0,\infty)$ and
$r\geq1$, there exists a constant\break $c_4(\beta^*, L^*,r) > 0$, such that
the new density estimator $\hat{p}_n$ with adaptively chosen bandwidth
according to \eqref{lepski_minimization} satisfies
\[
\sup_{(\beta,L)\in\mathcal{R}(\beta^*,L^*)} \sup_{p \in\mathscr{P}_d(\beta,L)} \sup_{t \in\mathbb
{R}^d}
\E_p^{\otimes n} \biggl( \frac{\vert\hat{p}_n(t) -
p(t) \vert}{\tilde{\psi}^n_{p(t),\beta}} \biggr)^r
\leq c_4 \bigl(\beta^*,L^*,r \bigr)
\]
for all $n \geq2$, where
\[
\tilde{\psi}^n_{x,\beta} := \bigl[n^{-{\barbeta}/{(\barbeta+1)}}\vee(x/n
)^{{\barbeta}/{(2\barbeta+1)}} \bigr] (\log n)^{3/2}.
\]
\end{thmm}

\begin{figure}[t]

\includegraphics{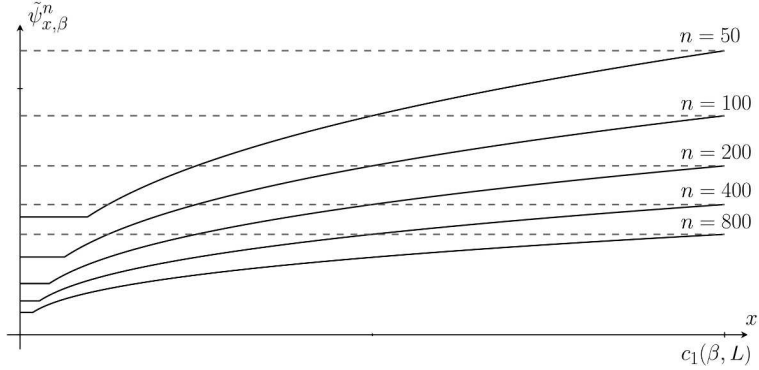}

\caption{New upper bound without logarithmic factor (solid line),
classical upper bound (dashed line).}
\label{fig:upperbound}
\end{figure}

The $p$-dependent\vspace*{1pt} speed of convergence $\tilde{\psi}_{p(t),\beta}^n$
(except the logarithmic factor) is plotted in Figure~\ref
{fig:upperbound}, which shows the superiority of the new estimator in
low density regions. It also depicts that the new estimator is able to
adapt to regime (ii) up to a logarithmic factor, and that it improves
the rate of convergence significantly in both regimes as compared to
the classical mini\-max rate. Besides, although not emphasized before,
$\hat{p}_n$ is fully adaptive to the smoothness in terms of H\"older
regularity.
%
%

As $\psi$ and $\tilde{\psi}$ coincide (up to a logarithmic factor) in
regime (ii) but differ in regime (i), the question arises whether the
breakpoint
\[
n^{-\barbeta/(\barbeta+1)}
\]
describes the fundamental bound on the range of adaptation to lowest
density regions. The following result shows that this is indeed the
case as long as the density's regularity is unknown.

%
\begin{thmm} \label{superefficiency}
For any $\beta_2 < \beta_1 \leq2$ and any sequence $(\rho(n))$
converging to infinity with
\[
\rho(n) = O \bigl( n^{\frac{\beta_1-\beta_2}{(2\beta_1+1)(\beta_2+1)}}
(\log n)^{-3/2} \bigr),
\]
there exist $L_1, L_2 > 0$ and densities $p_n \in\mathscr
{P}_1(\beta_1,L_1)$ with
\[
\frac{n^{-\beta_1/(\beta_1+1)}}{p_n(t)}=o(1)
\]
as $n\rightarrow\infty$, such that for every estimator $T_n(t)$ satisfying
%
%
\begin{equation}
\label{bigrate} \E_{p_n}^{\otimes n}\bigl \vert T_n(t) -
p_n(t) \bigr\vert\leq c_4 \bigl(\beta_1^*,L_1^*,r
\bigr) \biggl( \frac{p_n(t) }{n} \biggr)^{{\beta
_1}/{(2\beta
_1+1)}} (\log n)^{3/2},
\end{equation}
there exist $n_0(\beta_1,\beta_2, L_1, L_2)$ and a constant $c>0$ both
independent of $t$, with
\[
\mathop{\sup_{q \in\mathscr{P}_d(\beta_2,L_2):}}_{ q(t) \leq c(n)
\cdot n^{-{\beta_2}/{(\beta_2+1)}}} \frac{\E_{q}^{\otimes n}
\vert T_n(t) - q(t) \vert}{n^{-{\beta_2}/{(\beta_2+1)}}} \geq c
\]
for all $n \geq n_0(\beta_1, \beta_2, L_1, L_2)$ and any sequence
$(c(n))$ with $c(n) \geq\rho(n)^{-1}$.
\end{thmm}

The following consideration provides a heuristic reason why adaptation
to regime (i) is not possible in principle. Consider the univariate and
Lipschitz continuous triangular density $p: \mathbb{R} \rightarrow
\mathbb{R}, x \mapsto(1-\vert x\vert) \mathbh{1}\{\vert x \vert
\leq
1 \}$. If $\delta_n < n^{-\beta/ (\beta+1)} = n^{-1/2}$, the expected
number of observations in $ \{ p \leq\delta_n \}$ is less
than one. Without the knowledge of the regularity, it is intuitively
clear that it is impossible to predict whether local averaging is
preferable to just estimating by zero.

\subsubsection{Adaptation to lowest density regions when \texorpdfstring{$\beta$}{beta} is
known} \label{ext1} If the H\"older exponent $\beta\in(0,2]$ is known
to the statistician, the form of the oracle estimator~\eqref{oracle}
suggests that some further improvement in regime (i) might be possible
by considering the truncated estimator
%
%
\begin{equation}
\label{newest} \hat{p}_n(\cdot) \cdot\mathbh{1} \bigl\{
\hat{p}_n(\cdot) \geq n^{-{\barbeta}/{(\barbeta+1)}} (\log n)^{\zeta
_1} \bigr\}
\end{equation}
for some suitable constant $\zeta_1 > 0$. In fact, elementary algebra
shows that this threshold does not affect the performance in regime
(ii) (up to a logarithmic term). For isotropic H\"older smoothness, we
prove in the supplemental article [\citet{supp}] that
the estimator \eqref{newest} indeed attains the $p$-dependent speed of
convergence
\[
\vartheta_{p(t),\beta}^n = \psi_{p(t),\beta}^n
\vee n^{-\zeta_2}
\]
up to logarithmic terms, with $ \psi_{x,\beta}^n$ as defined in
Theorem~\ref{lowerbounddensity}. Here, the constant $\zeta_2$ can be made
arbitrarily large by enlarging $c_3$ and $\zeta_1$. That is, if the H\"
older exponent is known, adaptation to regime (i) is possible to a
large extent.

\subsubsection{Extension to \texorpdfstring{$\beta>2$}{beta>2}} \label{ext2} As concerns an
extension of Theorems \ref{lowerbounddensity} and \ref
{upperbounddensity} to arbitrary $\beta>2$, Lemma~\ref{varapprox}(ii)
demonstrates that the variance of the kernel density estimator never
falls below the reference speed of convergence $\tilde{\psi
}_{p(t),\beta
}^n$. However, it can be substantially larger, resulting in a lower
speed of convergence as compared to the reference speed of convergence.
Therefore, it seems necessary to introduce a $p$-dependent speed of
convergence which does not incorporate the value of the density $p(t)$
only but also information on the derivatives. An exception of
outstanding importance are points $t$ close to the support boundary,
because not only $p(t)$ itself but also all derivatives are necessarily
small. Theorem A.1, which is deferred to the supplemental article
[\citet{supp}], reveals that our procedure then even
reaches the fast adaptive speed of convergence at the support boundary
for \textit{every} $\beta>0$. In fact, as $\beta\rightarrow\infty$,
adaptive rates arbitrarily close to $n^{-1}$ can be attained.

\section{Application to support recovery} \label{sec:supportapplication}

The phenomenon of faster rates of convergence in regions where the
density is small may have strong consequences on plug-in rules for
certain functionals of the density. As an application of the results of
Section~\ref{sec:densitybounds}, we investigate the support
plug-in functional. Support estimation has a long history in the
statistical literature. \citet{geffroy1964} and
R\'enyi and Sulanke (\citeyear{renyisulanke1963},
\citeyear{renyisulanke1964}) are cited as pioneering reference most
commonly, followed by
further contributions of \citet{chevalier1976}, \citet
{devroyewise1980},
\citet{grenander1981}, \citet{hall1982}, \citet{groeneboom1988},
Tsybakov (\citeyear{tsybakov1989}, \citeyear{tsybakov1991},
\citeyear{tsybakov1997}), \citet{cuevas1990}, Korostelev and
Tsybakov (\citeyear{korostelevtsybakov1993}), \citet
{haerdleparktsybakov1995}, \citet{mammentsybakov1995}, \citet
{cuevasfraiman1997}, \citet{gayraud1997}, \citet
{hallnussbaumstern1997}, Ba\'illo,
Cuevas and Justel (\citeyear{baillocuevasjustel2000}), Cuevas and Rodr\'iguez-Casal
(\citeyear{cuevasrodriguez-casal2004}), Klemel\"
a (\citeyear{klemela2004}), and \citet{biaucadrepelletier2008},
Biau, Cadre, Mason and
Pelletier (\citeyear{biaucadremasonpelletier2009}), \citet
{brunel2013} and \citet{cholaquidiscuevasfraiman2014} as a by far
nonexhaustive list of contributions.
In order to demonstrate the substantial improvement in the rates of
convergence for the plug-in support estimator based on the new density
estimator, we first establish minimax lower bounds for support
estimation under the margin condition which have not been provided in
the literature so far. Theorems \ref{lowerboundsupport} and
\ref{upperboundsupport} then reveal that the minimax rates
for the support estimation problem are substantially faster than for
the level set estimation problem, as already conjectured in \citet
{tsybakov1997}. In fact, in the level set estimation framework, when
$\beta$
and $L$ are given, the classical choice of a bandwidth of order
$n^{-1/(2\beta+d)}$ in case of isotropic H\"older smoothness leads
directly to a minimax-optimal plug-in level set estimator as long as
the offset is suitably chosen [\citet{rigolletvert2009}]. In contrast,
this bandwidth produces suboptimal rates in the support estimation
problem, no matter how the offset is chosen. At first sight, this makes
the plug-in rule as a by-product of density estimation inappropriate.
We shall demonstrate subsequently, however, that our new density
estimator avoids this problem.
In order to line up with the results of \citet{cuevasfraiman1997} and
\citet{rigolletvert2009}, we work essentially under the same type of
conditions. The distance between two subsets $A$ and $B$ of $\mathbb
{R}^d$ is measured by
\[
d_\Delta(A,B) := \bblambda^d (A \Delta B),
\]
where $\Delta$ denotes the symmetric difference of sets
\[
A \Delta B := (A \setminus B) \cup(B \setminus A).
\]
Subsequently, $\bar{A}$ denotes the topological closure of a set
$A\subset\R^d$. We impose the following condition, which
characterizes the complexity of the problem.
It was introduced by \citet{polonik1995} [see also \citet
{mammentsybakov1999}, \citet{tsybakov2004} and \citet
{cuevasfraiman1997}], where the
latter authors referred to it as sharpness order.

%
\begin{definition}[(Margin condition)] \label{margincondition}
A density $p: \mathbb{R}^d \rightarrow\mathbb{R}$ is said to satisfy
the $\kappa$-margin condition with exponent $\gamma> 0$, if
\[
\bblambda^d \bigl( \overline{ \bigl\{ x \in\mathbb{R}^d
\vert0 < p(x) \leq\varepsilon\bigr\}} \bigr) \leq\kappa_2 \cdot
\varepsilon^\gamma
\]
for all $0 < \varepsilon\leq\kappa_1$, where $\kappa=
(\kappa
_1, \kappa_2) \in(0,\infty)^2$.
\end{definition}

In particular, $\bblambda^d(\partial\Gamma_p) = 0$ for every density
which satisfies the margin condition, where $\partial\Gamma_p$ denotes
the boundary of the support $\Gamma_p$.
To highlight the line of ideas, we restrict the application to the
important special case of isotropic smoothness. Let $\mathscr
{H}_d^{\mathrm{iso}}(\beta,L)$ denote the isotropic H\"older class with
one-dimensional parameters $\beta$ and $L$, which is for $0 < \beta
\leq1$ defined by
\[
\mathscr{H}_d^{\mathrm{iso}}(\beta,L) := \bigl\{ f:
\mathbb{R}^d \rightarrow\mathbb{R} : \bigl\vert f(x) - f(y) \bigr\vert\leq L
\Vert x-y \Vert_2^\beta\mbox{ for all } x,y\in
\R^d \bigr\}.
\]
For $\beta> 1$, it is defined as the set of all functions $f: \mathbb
{R}^d \rightarrow\mathbb{R}$ that are $\lfloor\beta\rfloor$ times
continuously differentiable such that the following property is satisfied:
%
%
\begin{equation}
\label{bedingung} \bigl\llvert f(x) - P_{y,\lfloor\beta\rfloor
}^{(f)}(x) \bigr
\rrvert\leq L \Vert x-y \Vert_2^\beta\qquad\mbox{for all } x,y
\in\mathbb{R}^d,
\end{equation}
where
\[
P_{y,l}^{(f)}(x) := \sum_{\vert k \vert\leq l}
\frac{D^k f(y)}{k_1!
\cdots k_d!} (x_1-y_1)^{k_1}
\cdots(x_d-y_d)^{k_d}
\]
with $\vert k \vert:= \sum_{i=1}^d k_i$ and the partial differential operator
\[
D^k := \frac{\partial^{\vert k \vert} }{\partial x_1^{k_1} \cdots
\partial x_d^{k_d}}
\]
denotes the multivariate Taylor polynomial of $f$ at the point $y \in
\mathbb{R}^d$ up to the $l$th order; see also \eqref
{taylorpolynomial} for the coinciding definition in one dimension.
Correspondingly, $\mathscr{P}_d^{\mathrm{iso}}(\beta,L)$ denotes the set of
probability densities contained in $\mathscr{H}_d^{\mathrm{iso}}(\beta,L)$. The
following lemma demonstrates that not every combination of margin
exponent and H\"older continuity is possible.

%
\begin{lemma} \label{ae}
There exists a compactly supported density in $\mathscr
{P}_d^{\mathrm{iso}}(\beta,L)$ satisfying a margin condition to the exponent
$\gamma$ if and only if $\gamma\beta\leq1$.
\end{lemma}

\subsection{Lower risk bounds for support recovery} \label
{subsec:supportbounds}

For any subset $A \subset\mathbb{R}^d$ and $\varepsilon> 0$, the
closed outer parallel set of $A$ at distance $\varepsilon>0$ is given by
\[
A^\varepsilon:= \Bigl\{ x \in\mathbb{R}^d :
\inf_{y \in A} \Vert x - y \Vert_2 \leq\varepsilon\Bigr\}
\]
and the closed inner $\varepsilon$-parallel set by $
A^{-\varepsilon} := \overline{ ( ( A^c
)^\varepsilon
)^c}$. Here, $\Vert\cdot\Vert_2$ denotes the Euclidean norm (on
$\R^d$). A support satisfying
\[
0 < \liminf_{\varepsilon\rightarrow0} \frac
{\bblambda^d (\Gamma_p \setminus\Gamma_p^{-\varepsilon
})}{\bblambda^d
(\Gamma_p^\varepsilon\setminus\Gamma_p)} \leq\limsup_
{\varepsilon\rightarrow0}
\frac{\bblambda^d (\Gamma_p
\setminus\Gamma_p^{-\varepsilon})}{\bblambda^d (\Gamma
_p^\varepsilon
\setminus\Gamma_p)} < \infty
\]
is referred to as boundary regular support. Note that a support is
always boundary regular if its Minkowski surface measure is
well-defined (in the sense that outer and inner Minkowski content exist
and coincide). The minimax lower bound is formulated under the
assumption of $\Gamma_p$ fulfilling the following complexity condition
(to the exponent $\mu= \gamma\beta$), which even slightly weakens the
assumption of boundary regularity under the margin condition.

%
\begin{definition}[(Complexity condition)] \label{cc}
A set $A$ is said to satisfy the $\xi$-complexity condition to the
exponent $\mu>0$ if for all $0 < \varepsilon\leq\xi_1$ there exists
a disjoint decomposition $A = A_{1,\varepsilon} \cup A_{2,\varepsilon}$
such that
\[
\frac{\bblambda^d(A_{1,\varepsilon}^\varepsilon\setminus
A_{1,\varepsilon}) \vee\bblambda^d(A_{2,\varepsilon
})}{\varepsilon^{\mu}} \leq\xi_2,
\]
where $\xi= (\xi_1, \xi_2) \in(0,\infty)^2$.
\end{definition}

Note that a boundary regular support of a $(\beta,L)$-H\"older-smooth
density satisfying the margin condition to the exponent $\gamma$
fulfills the complexity condition to the exponent $\mu\geq\gamma
\beta
$ for the canonical decomposition $\Gamma_p = \Gamma_p \cup\varnothing$.
Let us finally relate the margin condition (\ref{margincondition}) to the
two-sided margin condition
\[
\bblambda^d \bigl\{ x \in\mathbb{R}^d : 0 < \bigl\vert p(x) -
\lambda\bigr\vert\leq\varepsilon\bigr\} \leq c \varepsilon^\gamma,
\]
which is imposed in the context of density level set estimation for
some level $\lambda> 0$; cf. \citet{rigolletvert2009}. If $\Gamma
_{p,\lambda}= \{ x \in\R^d : p(x) > \lambda\}$ denotes the $\lambda
$-level set at level $\lambda>0$, the two-sided $(\kappa,\gamma
)$-margin condition provides the bound
%
%
\begin{equation}
\label{endspurt!} \bblambda^d \bigl( \Gamma_{p,\lambda}^\varepsilon
\setminus\Gamma_{p,\lambda} \bigr) \leq\kappa_2 \bigl(c
\varepsilon^{\beta\wedge
1} \bigr)^\gamma
\end{equation}
for all $\varepsilon\leq\kappa_1$, where $c=L$ for $\beta\leq1$ and
$c= \sup_{x \in\R^d} \Arrowvert\nabla p (x) \Arrowvert_2$ for
$\beta
>1$. In contrast, the margin condition at $\lambda=0$ provides no bound
on $\bblambda^d(\Gamma_p^\varepsilon\setminus\Gamma_p)$. The
complexity condition is a mild assumption which guarantees such type of
bound. For $\beta\leq1$, the relation \eqref{endspurt!} for $\lambda
=0$ implies the complexity condition to the exponent $\mu= \gamma
\beta
$. Note that the typical situation is indeed
\[
\bblambda^d \bigl(\Gamma_p^\varepsilon\setminus
\Gamma_p \bigr) / \varepsilon= \mathcal{O}(1)\quad \mbox{and}\quad
\varepsilon/\bblambda^d \bigl(\Gamma_p^\varepsilon
\setminus\Gamma_p \bigr) = \mathcal{O}(1)
\]
as $\varepsilon\rightarrow0$. For instance, this holds true for any
finite union of compact convex sets in $\mathbb{R}^d$ as a consequence
of the isoperimetric inequality [Theorem III.2.2, \citet{chavel2001}] and
Theorem~3.1 [Bhattacharya and Rango Rao (\citeyear
{bhattacharyarangorao1976})]. If it exists, the limit
\[
\lim_{\varepsilon\searrow0} \frac{\bblambda^d(\Gamma
_p^\varepsilon
\setminus\Gamma_p)}{\varepsilon}
\]
corresponds to the surface measure of the boundary if the latter is
sufficiently regular. Due to the relation $\gamma\beta\leq1$ by
Lemma~\ref{ae} and the decomposition into suitable subsets, the
complexity condition relaxes this regularity condition on the surface
area substantially. The subset of $\mathscr{P}_d^{\mathrm{iso}}(\beta, L)$
consisting of densities satisfying the $\kappa$-margin condition to the
exponent $\gamma$ with support fulfilling the $\xi$-complexity
condition to the exponent $\mu= \gamma\beta$ is denoted by $\mathscr
{P}_d^{\mathrm{iso}}(\beta, L, \gamma, \kappa,\xi)$.

%
\begin{thmm}[(Minimax lower bound)] \label{lowerboundsupport}
For any $\beta> 0$ and any margin exponent $\gamma> 0$ with $\gamma
\beta\leq1$, there exist $c_5(\beta,L)>0$, $n_0(\beta,L,\gamma)
\in
\mathbb{N}$ and parameters $\kappa, \xi\in(0,\infty)$, such that the
minimax risk with respect to the measure of symmetric difference of
sets is bounded from below by
\[
\inf_{\hat{\Gamma}_n} \sup_{p \in\mathscr
{P}_d^{\mathrm{iso}}(\beta, L, \gamma, \kappa, \xi)} \E_p^{\otimes n}
\bigl[ d_\Delta(\hat{\Gamma}_n, \Gamma_p)
\bigr] \geq c_5(\beta,L) \cdot n^{-{\gamma\beta}/{(\beta+d)}}
\]
for all $n \geq n_0(\beta,L,\gamma)$.
\end{thmm}

\subsection{Minimax-optimal plug-in rule}
We use the plug-in support estimator with the kernel density estimator
of Section~\ref{sec:densitybounds}. This density estimator improves the
rate of convergence
in particular at the support boundary. For the isotropic procedure, the
index set $\mathcal{J}$ is restricted to bandwidths coinciding in all
components. We even simplify the ordering by estimated variances in
condition \eqref{truncated_variance_estimator}
''for all $m \in\mathcal{J}$ with $\hat{\sigma}_t^2(m) \geq\hat
{\sigma
}_t^2(j)$'' by the classical order ``for all $m \in\mathcal{J} \mbox{
with } m\geq j$'' as Lemma~\ref{monotoneconvolution} shows that the
relevant orderings are equivalent up to multiplicative constants for $0
< \beta\leq2$. Furthermore, under isotropic smoothness it is natural
to use a rotation invariant kernel, that is, $K(x) = \tilde{K}(\Vert x
\Vert_2)$ with $\tilde{K}$ supported on $[0,1]$ and continuous on its
support with $\tilde{K}(0)>0$. The following theorem shows that the
corresponding plug-in rule
\[
\hat{\Gamma}_n = \overline{ \bigl\{ x \in\mathbb{R}^d :
\hat{p}_n(x) > \alpha_n \bigr\} }
\]
with offset level
%
%
\begin{equation}
\label{offset} \alpha_n := c_6(\beta,L) \biggl(
\frac{(\log n)^{3/2}}{n} \biggr)^{{\beta}/ {(\beta+d)}} \sqrt{\log n}
\end{equation}
and constant $c_6(\beta,L)$ specified in the proof of the following theorem,
is able to recover the support with minimax optimal rate, up to a
logarithmic factor.

%
\begin{thmm}[(Uniform upper bound)] \label{upperboundsupport}
For any $\beta\leq2$, $\gamma> 0$ with $\gamma\beta\leq1$ and
$\kappa, \xi\in(0,\infty)$, there exist a constant $c_7 = c_7(\beta
,L,\gamma,\kappa,\xi) > 0$ and $n_0 \in\mathbb{N}$, such that
\[
\sup_{p \in\mathscr{P}_d^{\mathrm{iso}}(\beta,L,\gamma,\kappa,\xi
)} \E_p^{\otimes n} \bigl[
d_\Delta( \hat{\Gamma}_n, \Gamma_p ) \bigr]
\leq c_7 \cdot n^{-{\gamma\beta}/{(\beta+d)}} (\log n)^{2\gamma}
\]
for all $n \geq n_0$.
\end{thmm}

As the rate already indicates, it is getting apparent from the proof
that this result can be established only if the minimax optimal density
estimator actually adapts up to the fastest rate in regime (ii).

%
\begin{remark}
The results show the simultaneous optimality of the adaptive density
estimator of Section~\ref{sec:densitybounds} in the plug-in rule for
support estimation. Correspondingly, they are restricted to $\beta\leq
2$. Whether the rate
$n^{-\gamma\beta/(\beta+d)}$ is minimax optimal for $\beta>2$ provided
$\gamma\beta\leq1$, and whether it can be attained by a plug-in rule
in principle, remains open for the moment.
\end{remark}

Let us finally point out two consequences. We have shown that the
optimal minimax rates for support estimation are significantly faster
than the corresponding rates for level set estimation
\[
n^{-{\gamma\beta}/{(2\beta+d)}}
\]
under the margin condition [\citet{rigolletvert2009}]. Although any level
set of a fixed density satisfying the margin condition to the exponent
$\gamma$ fulfills the complexity condition to the exponent $\mu=
\gamma\beta$ as long as $\beta\leq1$, the hypotheses in the proof of
the lower bounds of \citet{rigolletvert2009} do even satisfy this
condition for some fixed $\xi$, uniformly in $n$, as well. Hence, their
optimal minimax rates of convergence remain the same under our
condition. On an intuitive level, this phenomenon can be nicely
motivated by comparing the Hellinger distance $H(\mathbb{P},\mathbb
{Q})$ between the probability measure $\mathbb{P}$ with Lebesgue
density $p$ and $\mathbb{Q}$ whose Lebesgue density $q = p + \tilde{p}$
is a perturbation of $p$ with a small function $\tilde{p}$ around the
level $\alpha\geq0$; see \citet{tsybakov1997}, Extension (E4). If
$\alpha
> 0$, then simple Taylor expansion of $\sqrt{p + \tilde{p}}$ yields
$H^2(\mathbb{P},\mathbb{Q}) \sim\int\tilde{p}^2 \,d \bblambda^d$,
whereas $H^2(\mathbb{P},\mathbb{Q}) \sim\int\tilde{p} \,d
\bblambda^d$ in case $\alpha= 0$. Thus, perturbations at the boundary
($\alpha=0$) can be detected with the higher accuracy resulting in
faster attainable rates for support estimation than for level set
estimation. Moreover, the rates for plug-in support estimators already
established in the literature by \citet{cuevasfraiman1997} turn
out to
be always suboptimal in case of H\"older continuous densities of
boundary regular support. To be precise, \citet{cuevasfraiman1997}
establish in Theorem~1(c) a convergence rate under the margin
condition given in terms of $\rho_n=n^{\rho}$ and the offset level
$\alpha_n=n^{-\alpha}$ (in their notation), which are assumed to
satisfy $0<\alpha<\rho$ and their condition (R2), namely
\[
\rho_n\int\vert\hat{p}_n-p\vert\,\Leb=
o_{\P} (1) \quad\mbox{and}\quad \rho_n\alpha_n^{1+\gamma}
= o(1) \qquad\mbox{as } n\rightarrow\infty.
\]
As a consequence, $\rho_n=o(n^{\beta/(2\beta+d)})$ for typical
candidates $p\in\mathscr{P}_d^{\mathrm{iso}}(\beta,L)$, that is, densities $p$
which are locally not smoother than $(\beta,L)$-regular. Under the
margin condition to the exponent $\gamma>0$, this limits their rate of
convergence $n^{-\rho+\alpha}$ to
\[
d_{\Delta}(\Gamma_p,\hat{\Gamma}_n) =
o_{\P} \bigl(n^{-({\beta}/{(2\beta+d)})({\gamma}/{(1+\gamma)})} \bigr),
\]
which is substantially slower than the above established minimax rate.
The crucial point is that even with the improved density estimator of
Section~\ref{sec:densitybounds}, the above mentioned condition on
$\rho
_n$ in (R2) cannot be improved, because any estimator can possess the
improved performance at lowest density regions only. For this reason,
the $L_1$-speed of convergence of a density estimator is not an
adequate quantity to characterize the performance of the corresponding
plug-in support estimator.

\section{Lemmas \texorpdfstring{\protect\ref{varapprox}}{5.1}--\texorpdfstring{\protect\ref{biasdiff}}{5.7}, proofs
of Theorems \texorpdfstring{\protect\ref{upperbounddensity}}{3.3} and \texorpdfstring{\protect\ref{superefficiency}}{3.4}} \label
{sec:proofs}

Due to space constraints, all remaining proofs are deferred to the
supplemental article\break [Patschkowski and Rohde (\citeyear{supp})]. In the proof of
Theorem~\ref{upperbounddensity}, we frequently make use of the
bandwidth
%
%
\begin{equation}
\label{new_opt_bw} \bar{h}_i := c_{8}(\beta,L) \cdot\max
\biggl\{ \biggl( \frac{\log n}{n} \biggr)^{({\barbeta}/{(\barbeta+1)})
({1}/{\beta_i})}, \biggl(
\frac
{p(t) \log n}{n} \biggr)^{({\barbeta}/{(2 \barbeta+ 1)})
({1}/{\beta
_i})} \biggr\}\hspace*{-30pt}
\end{equation}
for $i=1,\ldots,d$, with constant $c_{8}(\beta,L)$ of Lemma~\ref
{varapprox}, which can be thought of as an optimal adaptive bandwidth.
The truncation in the definition of $\bar{h}$ results from the
necessary truncation in $\sigma_{t,\mathrm{trunc}}^2$.
With the exponents
%
%
\begin{equation}
\label{grid_approx} \bar{j}_i = \bar{j}_i(t) := \biggl
\lfloor\log_2 \biggl( \frac
{1}{\bar
{h}_i} \biggr) \biggr\rfloor+ 1, \qquad i =
1, \ldots, n
\end{equation}
the bandwidth $ 2^{-\bar{j}_i}$ is an approximation of $\bar{h}_i$ by
the next smaller bandwidth on the grid $\mathcal{G}$ such that $\bar
{h}_i /2 \leq2^{-\bar{j}_i} \leq\bar{h}_i$ for all $i=1,\ldots,d$.

Before turning to the proof of Theorem~\ref
{upperbounddensity}, we collect some technical ingredients. First,
recall the classical upper bound on the bias of a kernel density
estimator. With the notation provided in Section~\ref
{sec:preliminairies}, and $K$ of order $\max_i \beta_i$ at least, we obtain
\begin{eqnarray*}
b_t(h)& :=& p(t) - \E_p^{\otimes n}
\hat{p}_{n,h}(t)\\
 &= &\int K(x) \bigl( p(t + h x) - p(t) \bigr) \,\Leb(x)
\\
&= &\sum_{i=1}^d \int K(x) \bigl( p
\bigl([t, t + h x]_{i-1} \bigr) - p \bigl([t, t + h x]_i
\bigr) \bigr) \,\Leb(x),
\end{eqnarray*}
using the notation $[x,y]_0 = y$, $[x,y]_d = x$, $[x,y]_i = (x_1,
\ldots
, x_i, y_{i+1}, \ldots, y_d)$, $i = 1, \ldots, d-1$ for two
vectors $x,y \in\mathbb{R}^d$ and denoting by $hx = (h_1 x_1, \ldots,
h_d x_d)$ the componentwise product.
Taylor expansions for those components $i$ with $\beta_i \geq1$ lead to
\begin{eqnarray*}
&&p \bigl([t,t+hx]_{i-1} \bigr) - p \bigl([t,t+hx]_i
\bigr)\\
&&\qquad= \sum_{k=1}^{\lfloor\beta_i
\rfloor}
p_{i,[t, t + h x]_i}^{(k)}(t_i) \frac{(h_i x_i)^k}{k!}
\\
&&\qquad\quad{} + \bigl( p \bigl( [t,t+hx]_{i-1} \bigr) - P_{t_i, \lfloor\beta_i
\rfloor
}^{(p_{i,[t,t+hx]_i})}
(t_i + h_i x_i) \bigr).
\end{eqnarray*}
Hence,
%
%
\begin{equation}
\label{bias} \bigl\vert b_t(h) \bigr\vert\leq L \sum
_{i=1}^d c_{9,i}(\beta)
h_i^{\beta_i} =: B_t(h)
\end{equation}
with constants
$c_{9,i}(\beta) := \int\vert x_i \vert^{\beta_i} \vert K(x) \vert
\,\Leb(x) < \infty$.

With a slight abuse of notation, dependencies on some
bandwidth $h = 2^{-j}$ are subsequently expressed in terms of the
corresponding grid exponent $j = (j_1, \ldots, j_d)$, that is, $B_t(h)$
equals $B_t(j)$, etc. For any multi-index $j$, we use the abbreviation
\[
\vert j \vert:= \sum_{i=1}^d
j_i.
\]
The following lemmata are crucial ingredients for the proof of
Theorem~\ref{upperbounddensity}.

%
\begin{lemma} \label{varapprox} \textup{(i)} For any $(\beta,L)$ with $0<\beta_i
\leq2$, $p \in\mathscr{P}_d(\beta,L)$, and for any bandwidth
$h=(h_1,\ldots,h_d)$ with $h_i \leq c_{8}(\beta,L) p(t)^{1/\beta
_i}$, $i=1,\ldots,d$ with
\[
c_{8}(\beta,L) := \min_{i=1,\ldots,d} \biggl(
\frac{2dL}{\Vert K
\Vert
_2^2} \int\vert x_i \vert^{\beta_i}
K^2(x) \,\Leb(x) \biggr)^{-1/\beta_i},
\]
the following inequality chain holds true
\[
\frac{1}{2} \frac{\Vert K \Vert_2^2}{n \prod_{i=1}^d h_i} p(t) \leq
\frac{1}{n}
\bigl((K_h)^2 \ast p \bigr) (t) \leq\frac{3}{2}
\frac
{\Vert K \Vert_2^2}{n \prod_{i=1}^d h_i} p(t).
\]

\textup{(ii)} For any constant $c_{10}>0$, there exists a
constant $c_{11}(\beta,L) =\break c_{11}(\beta,L, c_{10}) >0$, such that for
any $(\beta,L)$, $0<\beta_i < \infty$, $i=1,\ldots,d$, and $p \in
\mathscr{P}_d(\beta,L)$,
\[
\frac{c_{11}(\beta,L)}{n \prod_{i=1}^d h_i} p(t) \leq\frac
{1}{n} \bigl((K_h)^2
\ast p \bigr) (t)
\]
for every bandwidth $h=(h_1,\ldots, h_d)$ with $h_i \leq c_{10}
p(t)^{1/\beta_i}$, $i=1,\ldots, d$.

\textup{(iii)} For any density $p$ with isotropic H\"older smoothness $(\beta
,L)$, $0<\beta< \infty$ and bandwidth $h$, we have
\[
\frac{1}{n} \bigl((K_h)^2 \ast p \bigr) (t) \leq
\frac{L \Vert K \Vert_2^2}{n
h^d} \Bigl( h + \inf_{y \in\Gamma_p^c} \Vert t-y
\Vert_2 \Bigr)^{\beta},
\]
where $K$ is a rotation invariant kernel supported on the closed
Euclidean unit ball.
\end{lemma}

Lemma \ref{varapprox}(ii) provides an extension of the results of
Rohde (\citeyear{rohde2008,rohde2011}).

%
\begin{lemma}\label{monotoneconvolution}
There exists some constant $c_{12}(\beta,L)>0$, such that for any
$p\in
\mathscr{P}_d(\beta,L)$, $0 < \beta_i \leq2$, $i=1,\ldots,d$, and
$t\in
\R^d$ the inequality
\[
\sigma^2_{t,\mathrm{trunc}}(j \wedge m) \leq c_{12}(\beta,L)
\bigl( \sigma^2_{t,\mathrm{trunc}}(j) \vee\sigma^2_{t,\mathrm{trunc}}(m)
\bigr)
\]
holds true for all (nonrandom) indices $j = (j_{1}, \ldots, j_{d})$
and $m = (m_1, \ldots, m_d)$ with $j \geq\bar{j}$ componentwise. If
additionally $m \geq j$ componentwise, then
\[
\sigma_{t,\mathrm{trunc}}^2(j) \leq c_{12}(\beta,L) \sigma
_{t,\mathrm
{trunc}}^2(m).
\]
\end{lemma}

The next lemma carefully analyzes the ratio of the truncated quantities
$\sigma_{t,\mathrm{trunc}}^2$ and $\tilde{\sigma}_{t,\mathrm{trunc}}^2$.

%
\begin{lemma} \label{sigmabreve_to_sigmacheck}
For the quantities $\sigma_{t,\mathrm{trunc}}^2(h)$ and $\tilde
{\sigma
}_{t,\mathrm{trunc}}^2(h)$ defined in \eqref{sigma_trunc} and any
$\eta
\geq0$ holds
\[
\mathbb{P}^{\otimes n} \biggl( \biggl\llvert\frac{\tilde{\sigma
}_{t,\mathrm{trunc}}^2(h)}{\sigma_{t,\mathrm{trunc}}^2(h)} - 1 \biggr
\rrvert\geq\eta\biggr) \leq2 \exp\biggl( -\frac{3 \eta
^2}{2(3+2\eta) \Vert K \Vert_\mathrm{sup}^2}
\log^2 n \biggr).
\]
\end{lemma}

%
\begin{lemma} \label{lemma_hopt}
For any $(\beta,L)$ with $0<\beta_i \leq2$, $i=1,\ldots,d$, there
exist constants $c_{13}(\beta,L)$ and $c_{14}(\beta,L) > 0$ such that
for the multi-index $\bar{j}$ as defined in \eqref{grid_approx} and the
bias upper bound $B_t$ as given in \eqref{bias},
%
%
\begin{eqnarray}
B_t(\bar{j}) &\leq &c_{13}(\beta,L) \sqrt{
\sigma_{t,\mathrm
{trunc}}^2(\bar{j}) \log n}, \label{inequ1}
\\
\sqrt{\sigma_{t,\mathrm{trunc}}^2(\bar{j}) } &\leq&
c_{14}(\beta,L) \biggl\{ \biggl(\frac{\log n}{n}
\biggr)^{{\barbeta}/{(\barbeta+1)}} \vee\biggl( \frac{p(t) \log
n}{n} \biggr)^{{\barbeta}/{(2\barbeta+
1)}}
\biggr\}. \label{inequ2}
\end{eqnarray}
\end{lemma}

%
\begin{lemma} \label{tailsZ}
For any (nonrandom) index $j = (j_1, \ldots, j_d)$, the tail
probabilities of the random variable
\[
Y := \frac{\hat{p}_{n,j}(t) - \E_p^{\otimes n} \hat
{p}_{n,j}(t)}{\sqrt
{\sigma_{t,\mathrm{trunc}}^2(j) \log n}},
\]
are bounded by
\[
\mathbb{P}^{\otimes n} \bigl(\vert Y \vert\geq\eta\bigr) \leq2 \exp\biggl( -
\frac{\log n}{4} \cdot\bigl(\eta^2 \wedge\eta\bigr) \biggr)
\]
for any $\eta\geq0$, any $t \in\mathbb{R}^d$ and $n \geq n_0$ with
$n_0$ depending on $\Vert K \Vert_{\sup}$ only.
\end{lemma}

%
\begin{lemma} \label{tailsZ2}
Let $Z$ be some nonnegative random variable satisfying
\[
\mathbb{P} (Z \geq\eta) \leq2 \exp( -A \eta)
\]
for some $A>0$. Then
\[
\bigl( \E Z^m \bigr)^{1/m} \leq c_{15}
\frac{m}{A}
\]
for any $m \in\mathbb{N}$, where the constant $c_{15}$ does not depend
on $A$ and $m$.
\end{lemma}

%
\begin{lemma} [{[\citet{klutchnikoff2005}]}] \label{biasdiff}
For all $k,l \in\mathcal{J}$, the absolute value of the difference of
bias terms is bounded by
\[
\bigl\vert b_t(k \wedge l) - b_t(l) \bigr\vert\leq2
B_t(k)
\]
for all $t \in\mathbb{R}^d$.
\end{lemma}

\begin{pf*}{Proof of Theorem~\ref{upperbounddensity}}
Recall the notation of Section~\ref{sec:densitybounds} and denote
$\hat
{p}_{n,\hat{j}} = \hat{p}_n$. In a first step, the risk
\[
\E_p^{\otimes n} \bigl\vert\hat{p}_{n,\hat{j}}(t) - p(t)
\bigr\vert^r
\]
is decomposed as follows:
%
%
\begin{eqnarray}
\label{r+r-} &&\E_p^{\otimes n} \bigl\vert\hat{p}_{n,\hat{j}}(t)
- p(t) \bigr\vert^r \nonumber\\
&&\qquad= \E_p^{\otimes n} \bigl[ \bigl\vert
\hat{p}_{n,\hat{j}}(t) - p(t) \bigr\vert^r \cdot\mathbh{1} \bigl\{
\hat{\sigma}_t^2(\hat{j}) \leq\hat{\sigma
}_t^2(\bar{j}) \bigr\} \bigr]
\nonumber
\\[-8pt]
\\[-8pt]
\nonumber
&&\qquad\quad{} + \E_p^{\otimes n} \bigl[ \bigl\vert\hat{p}_{n,\hat{j}}(t) -
p(t) \bigr\vert^r \cdot\mathbh{1} \bigl\{ \hat{\sigma}_t^2(
\hat{j}) > \hat{\sigma}_t^2(\bar{j}) \bigr\} \bigr]
\\
&&\qquad=: R^{+} + R^{-}.\nonumber
\end{eqnarray}
We start with $R^+$, which is decomposed again as follows:
%
%
\begin{eqnarray}
\label{S_decomp} %
R^+ &\leq&3^{r-1} \bigl( \E_p^{\otimes n}
\bigl[\bigl \vert\hat{p}_{n,\hat{j}}(t) - \hat{p}_{n,\hat{j} \wedge\bar{j}}(t)
\bigr\vert^r \cdot\mathbh{1} \bigl\{ \hat{\sigma}_t^2(
\hat{j}) \leq\hat{\sigma}_t^2(\bar{j}) \bigr\} \bigr]
\nonumber
\\
&&{} + \E_p^{\otimes n} \bigl[\bigl \vert\hat{p}_{n,\hat{j} \wedge\bar
{j}}(t) -
\hat{p}_{n,\bar{j}}(t) \bigr\vert^r \cdot\mathbh{1} \bigl\{ \hat{
\sigma}_t^2(\hat{j}) \leq\hat{\sigma}_t^2(
\bar{j}) \bigr\} \bigr]
\nonumber
\\[-8pt]
\\[-8pt]
\nonumber
&&{} +\E_p^{\otimes n} \bigl[ \bigl\vert\hat{p}_{n,\bar{j}}(t) -
p(t)\bigr \vert^r \cdot\mathbh{1} \bigl\{ \hat{\sigma}_t^2(
\hat{j}) \leq\hat{\sigma}_t^2(\bar{j}) \bigr\} \bigr]
\bigr) %
\\
&=: &3^{r-1} (S_1 + S_2 + S_3),\nonumber
\end{eqnarray}
where we used the inequality $(x+y+z)^r \leq3^{r-1} (x^r+y^r+z^r)$ for
all $x,y,z \geq0$. This decomposition bears the advantage that only
kernel density estimators with well-ordered bandwidths are compared. We
focus on the estimation of $S_1, S_2$ and $S_3$ and start with $S_2$
using the selection scheme's construction. Clearly, $\hat{j} \in
\mathcal{A}$ as defined in \eqref{admissiblebandwidths}. As a
consequence, the following inequality holds true:
\begin{eqnarray*}
S_2 &\leq&c_3^r \E_p^{\otimes n}
\biggl[ \bigl( \hat{\sigma}_t^2(\bar{j}) \log n
\bigr)^{r/2} \cdot\mathbh{1} \biggl\{ \biggl\llvert\frac{\tilde{\sigma
}_{t,\mathrm{trunc}}^2(\bar{j})}{\sigma
_{t,\mathrm{trunc}}^2(\bar{j})} -
1 \biggr\rrvert< 1 \biggr\} \biggr]
\\
&&{}+ c_3^r \E_p^{\otimes n} \biggl[
\bigl( \hat{\sigma}_t^2(\bar{j}) \log n
\bigr)^{r/2} \cdot\mathbh{1} \biggl\{ \biggl\llvert\frac
{\tilde
{\sigma}_{t,\mathrm{trunc}}^2(\bar{j})}{\sigma_{t,\mathrm
{trunc}}^2(\bar
{j})} -
1 \biggr\rrvert\geq1 \biggr\} \biggr]
\\
&\leq&2^{r/2} c_3^r \biggl( \min\biggl\{
\sigma_{t,\mathrm
{trunc}}^2(\bar{j}), \frac{\Vert K \Vert_2^2 c_1}{n 2^{-\vert\bar
{j} \vert}} \biggr\} \log n
\biggr)^{r/2}
\\
&&{}+ c_3^r \biggl( \frac{\Vert K \Vert_2^2 c_1}{n 2^{-\vert\bar{j}
\vert}} \log n
\biggr)^{r/2} \mathbb{P}^{\otimes n} \biggl( \biggl\llvert
\frac{\tilde{\sigma}_{t,\mathrm{trunc}}^2(\bar
{j})}{\sigma
_{t,\mathrm{trunc}}^2(\bar{j})} - 1 \biggr\rrvert\geq1 \biggr),
\end{eqnarray*}
where we used the condition in the indicator function in the first
summand to bound the estimated truncated variance $\tilde{\sigma
}_{t,\mathrm{trunc}}^2$ from above by $2 \sigma_{t,\mathrm{trunc}}^2$,
and additionally the upper truncation level in the second summand. By
the deviation inequality of Lemma~\ref
{sigmabreve_to_sigmacheck}, we can further estimate $S_2$ by
\begin{eqnarray*}
S_2 &\leq&2^{r/2} c_3^r \bigl(
\sigma_{t,\mathrm{trunc}}^2(\bar{j}) \log n \bigr)^{r/2}
\\
& &{}+ c_3^r \biggl( \frac{\Vert K \Vert_2^2 c_1}{n
2^{-\vert
\bar{j} \vert}} \log n
\biggr)^{r/2} \cdot2 \exp\biggl( -\frac{3}{10
\Vert K \Vert_\mathrm{sup}^2} \log^2
n \biggr).
\end{eqnarray*}
The second term is always of smaller order than the first term because
$2^{-\vert\bar{j} \vert} \leq1$ and, therefore, for $n \geq2$,
\[
\biggl( \frac{\Vert K \Vert_2^2 c_1}{n 2^{-\vert\bar{j} \vert}} \log n
\biggr)^{r/2} \cdot2 \exp\biggl( -
\frac{3}{10 \Vert K \Vert_\mathrm
{sup}^2} \log^2 n \biggr) \leq c \biggl( \frac{\log^3 n}{n^2 (
2^{-\vert\bar{j} \vert} )^2}
\biggr)^{r/2}
\]
for some constant $c$ depending on $c_1$, $r$ and the kernel $K$ only. Finally,
\[
S_2 \leq c(\beta,L) \bigl( \sigma_{t,\mathrm{trunc}}^2(
\bar{j}) \log n \bigr)^{r/2}.
\]
We will now turn to $S_3$, the third term in \eqref{S_decomp}. We split
the risk into bias and stochastic error. It holds
%
%
\begin{equation}
\label{decomp_s3} S_3 \leq\E_p^{\otimes n} \bigl(
\bigl\vert\hat{p}_{n,\bar{j}}(t) - \E_p^{\otimes n}
\hat{p}_{n,\bar{j}}(t) \bigr\vert+ B_t(\bar{j}) \bigr)^r
\end{equation}
and by Lemma~\ref{lemma_hopt}
%
%
\begin{equation}
\label{bias_jq} B_t(\bar{j}) \leq c_{13}(\beta,L) \sqrt
{\sigma_{t,\mathrm
{trunc}}^2(\bar{j}) \log n}.
\end{equation}
Denoting by
%
%
\begin{equation}
\label{defZ} Z_k := \frac{\hat{p}_{n,k}(t) - \E_p^{\otimes n} \hat
{p}_{n,k}(t)}{\sqrt{\sigma_{t,\mathrm{trunc}}^2(k) \log n}} \qquad\mbox{for
} k \in
\mathcal{J},
\end{equation}
the decomposition \eqref{decomp_s3}, the bias variance relation
\eqref{bias_jq} and the inequality $(x+y)^r \leq2^{r-1} (x^r+y^r)$,
$x,y \geq0$ together with Lemma~\ref{tailsZ2} yields
\begin{eqnarray*}
S_3 &\leq&\bigl( \sigma_{t,\mathrm{trunc}}^2(\bar{j}) \log
n \bigr)^{r/2} \cdot\E_p^{\otimes n} \bigl( \vert
Z_{\bar{j}} \vert+ c_{13}(\beta,L) \bigr)^r
\\
&\leq&\bigl( \sigma_{t,\mathrm{trunc}}^2(\bar{j}) \log n
\bigr)^{r/2} \cdot2^{r-1} \E_p^{\otimes n}
\bigl( \vert Z_{\bar{j}} \vert^r + c_{13}(
\beta,L)^r \bigr)
\\
&\leq& c(\beta,L) \bigl( \sigma_{t,\mathrm{trunc}}^2(\bar{j}) \log n
\bigr)^{r/2}.
\end{eqnarray*}
It remains to show an analogous result for $S_1$, the first term in
\eqref{S_decomp}. Clearly,
%
%
\begin{eqnarray}
\label{S1_inequ} %
S_1 &\leq&\sum
_{j \in\mathcal{J}} \E_p^{\otimes n} \bigl[ \bigl( \bigl\vert
\hat{p}_{n,j}(t) - \E_p^{\otimes n}
\hat{p}_{n,j}(t) \bigr\vert+ \bigl\vert\hat{p}_{n,j \wedge\bar{j}}(t) -
\E_p^{\otimes n} \hat{p}_{n,j \wedge\bar{j}}(t) \bigr\vert
\nonumber
\\[-8pt]
\\[-8pt]
\nonumber
&&{} + \bigl\vert b_t(j \wedge\bar{j}) - b_t(j) \bigr\vert
\bigr)^r \cdot\mathbh{1} \bigl\{ \hat{\sigma}_t^2(j)
\leq\hat{\sigma}_t^2(\bar{j}), \hat{j} = j \bigr\}
\bigr].
\end{eqnarray}
By Lemmas \ref{biasdiff} and \ref{lemma_hopt},
\[
\bigl\vert b_t(j \wedge\bar{j}) - b_t(j) \bigr\vert\leq2
B_t(\bar{j}) \leq2 c_{13}(\beta,L) \sqrt{
\sigma_{t,\mathrm{trunc}}^2(\bar{j}) \log n}.
\]
On account of this inequality and in view of \eqref{S1_inequ}, it
suffices to bound the expectations in the following expression:
%
%
\begin{eqnarray}
\label{stress} &&S_1 \leq 3^{r-1} \bigl(
\sigma_{t,\mathrm{trunc}}^2(\bar{j}) \log n \bigr)^{r/2}
\nonumber
\\
& &\qquad{}\times\biggl\{ \sum_{j \in\mathcal{J}} \E_p^{\otimes n}
\biggl[ \biggl( \frac{\vert\hat{p}_{n,j}(t) - \E_p^{\otimes n} \hat{p}_{n,j}(t)
\vert}{\sqrt{\sigma_{t,\mathrm{trunc}}^2(\bar{j}) \log n}} \biggr)^r
\mathbh{1} \bigl\{ \hat{
\sigma}_t^2(j) \leq\hat{\sigma}_t^2(
\bar{j}), \hat{j} = j \bigr\} \biggr]
\nonumber
\\[-8pt]
\\[-8pt]
\nonumber
&&\qquad{} + \sum_{j \in\mathcal{J}} \E_p^{\otimes n}
\biggl[ \biggl( \frac{\vert\hat{p}_{n,j \wedge\bar{j}}(t) - \E
_p^{\otimes n} \hat
{p}_{n,j \wedge\bar{j}}(t) \vert}{\sqrt{\sigma_{t,\mathrm
{trunc}}^2(\bar{j}) \log n}} \biggr)^r \mathbh{1} \bigl\{ \hat{
\sigma}_t^2(j) \leq\hat{\sigma}_t^2(
\bar{j}), \hat{j} = j \bigr\} \biggr]
\\
&&\qquad{} + \sum_{j \in\mathcal{J}} 2^r c_{13}(
\beta,L)^r \cdot\mathbb{P}^{\otimes n} ( \hat{j} = j) \biggr\}.
\nonumber
\end{eqnarray}
Denoting
%
%
\begin{equation}
\label{endspurt1} A_{j,\bar{j}} := \biggl\{ \biggl\llvert\frac{\tilde
{\sigma
}_{t,\mathrm
{trunc}}^2(j)}{\sigma_{t,\mathrm{trunc}}^2(j)} -
1 \biggr\rrvert< \frac
{1}{2} \mbox{ and } \biggl\llvert
\frac{\tilde{\sigma}_{t,\mathrm
{trunc}}^2(\bar{j})}{\sigma_{t,\mathrm{trunc}}^2(\bar{j})} - 1 \biggr
\rrvert< \frac{1}{2} \biggr\},
\end{equation}
it follows
\begin{eqnarray*}
&&\sum_{j \in\mathcal{J}} \E_p^{\otimes n}
\biggl[ \biggl( \frac
{\vert
\hat{p}_{n,j}(t) - \E_p^{\otimes n} \hat{p}_{n,j}(t) \vert}{\sqrt
{\sigma
_{t,\mathrm{trunc}}^2(\bar{j}) \log n}} \biggr)^r \mathbh{1} \bigl\{
\hat{
\sigma}_t^2(j) \leq\hat{\sigma}_t^2(
\bar{j}), \hat{j} = j \bigr\} \biggr]
\\
&&\qquad = \sum_{j \in\mathcal{J}} \E_p^{\otimes n}
\biggl[ \biggl( \frac{\vert\hat{p}_{n,j}(t) - \E_p^{\otimes n} \hat{p}_{n,j}(t)
\vert}{\sqrt{\sigma_{t,\mathrm{trunc}}^2(\bar{j}) \log n}} \biggr)^{r}
\mathbh{1} \bigl\{ \hat{
\sigma}_t^2(j) \leq\hat{\sigma}_t^2(
\bar{j}), \hat{j} = j \bigr\} \cdot\mathbh{1}_{A_{j,\bar{j}}} \biggr]
\\
&&\qquad\quad{} + \sum_{j \in\mathcal{J}} \E_p^{\otimes n}
\biggl[ \biggl( \frac{\vert\hat{p}_{n,j}(t) - \E_p^{\otimes n} \hat{p}_{n,j}(t)
\vert}{\sqrt{\sigma_{t,\mathrm{trunc}}^2(\bar{j}) \log n}} \biggr)^{r}
\mathbh{1} \bigl\{ \hat{
\sigma}_t^2(j) \leq\hat{\sigma}_t^2(
\bar{j}), \hat{j} = j \bigr\} \cdot\mathbh{1}_{A_{j,\bar{j}}^c} \biggr]
\\
&&\qquad =: S_{1,1} + S_{1,2}.
\end{eqnarray*}
Applying Lemma~\ref{tailsZ2} and H\"older's inequality for any $p>1$,
\begin{eqnarray*}
S_{1,1} &\leq&\biggl( \frac{3( 1 \vee c_1 \Vert K \Vert_2^2 )}{c_1
\Vert K \Vert_2^2} \biggr)^{r/2} \sum
_{j \in\mathcal{J}} \E_p^{\otimes
n} \bigl[
\vert Z_j \vert^r \cdot\mathbh{1} \{ \hat{j}=j \} \bigr]
\\
&\leq&\biggl( \frac{3( 1 \vee c_1 \Vert K \Vert_2^2 )}{c_1 \Vert K
\Vert_2^2} \biggr)^{r/2} \biggl( 1+ \sum
_{j \in\mathcal{J}} \E_p^{\otimes n} \bigl[ \vert
Z_j \vert^r \mathbh{1} \bigl\{ \vert Z_j \vert
\geq1 \bigr\} \mathbh{1} \{ \hat{j}=j \} \bigr] \biggr)
\\
&\leq&\biggl( \frac{3( 1 \vee c_1 \Vert K \Vert_2^2 )}{c_1 \Vert K
\Vert
_2^2} \biggr)^{r/2}\\
&&{}\times \biggl( 1+ \sum
_{j \in\mathcal{J}} \E_p^{\otimes n} \bigl[ \vert
Z_j \vert^{rp} \mathbh{1} \bigl\{ \vert Z_j
\vert\geq1 \bigr\} \bigr]^{1/p} \cdot\mathbb{P} ( \hat{j}=j )^{{(p-1)}/{p}}
\biggr)
\\
&\leq&\biggl( \frac{3( 1 \vee c_1 \Vert K \Vert_2^2 )}{c_1 \Vert K
\Vert
_2^2} \biggr)^{r/2} \biggl(
1+c_{15}^{r} \biggl( \frac{8rp}{\log
n}
\biggr)^r \sum_{j \in\mathcal{J}} \mathbb{P} (
\hat{j}=j )^{
{(p-1)}/{p}} \biggr)
\\
&\leq&\biggl( \frac{3( 1 \vee c_1 \Vert K \Vert_2^2 )}{c_1 \Vert K
\Vert
_2^2} \biggr)^{r/2} \biggl(
1+c_{15}^{r} \biggl( \frac{8rp}{\log
n}
\biggr)^r \biggl( \sum_{j \in\mathcal{J}} \mathbb{P} (
\hat{j}=j ) \biggr)^{{(p-1)}/{p}} \cdot\vert\mathcal{J} \vert^{{1}/{p}}
\biggr).
\end{eqnarray*}
By the constraint $2^{-\vert j \vert} \geq\log^2 n / n$ for any $j
\in
\mathcal{J}$, there exists some constant $c>0$ such that $\vert
\mathcal
{J} \vert\leq c (\log n)^d$. Setting finally $p = d \log n$, yields
$S_{1,1} \leq c(\beta^*,L^*)$. As concerns $S_{1,2}$, by the
Cauchy--Schwarz inequality,
\begin{eqnarray*}
S_{1,2} &\leq&\sum_{j \in\mathcal{J}} \biggl(
\frac{\sigma
_{t,\mathrm
{trunc}}^2(j)}{\sigma_{t,\mathrm{trunc}}^2(\bar{j})} \biggr)^{r/2} \E
_p^{\otimes n}
\bigl[ \vert Z_j \vert^{r} \mathbh{1} \{ \hat{j} = j \}
\mathbh{1}_{A_{j,\bar{j}}^c} \bigr]
\\
&\leq&\sum_{j \in\mathcal{J}} \biggl( \frac{\sigma_{t,\mathrm
{trunc}}^2(j)}{\sigma_{t,\mathrm{trunc}}^2(\bar{j})}
\biggr)^{r/2} \E_p^{\otimes n} \bigl[ \vert
Z_j \vert^{2r} \mathbh{1} \{ \hat{j}=j \}
\bigr]^{1/2}
\\
&&{}\times\biggl\{ \mathbb{P}^{\otimes n} \biggl( \biggl\llvert
\frac{\tilde{\sigma}_{t,\mathrm{trunc}}^2(j)}{\sigma_{t,\mathrm
{trunc}}^2(j)} - 1 \biggr\rrvert\geq\frac{1}{2} \biggr) + \mathbb
{P}^{\otimes n} \biggl( \biggl\llvert\frac{\tilde{\sigma}_{t,\mathrm
{trunc}}^2(\bar{j})}{\sigma_{t,\mathrm{trunc}}^2(\bar{j})} - 1 \biggr
\rrvert
\geq\frac{1}{2} \biggr) \biggr\}^{1/2}.
\end{eqnarray*}
Via the lower and upper truncation levels in the definition of $\sigma
_{t,\mathrm{trunc}}^2$,
%
%
\begin{equation}
\label{haesslich} \frac{\sigma_{t,\mathrm{trunc}}^2(k)}{\sigma
_{t,\mathrm{trunc}}^2(l)} \leq\frac{(1 \vee c_1 \Vert K \Vert_2^2)
n^2}{\log^4 n} \qquad\mbox{for any } k,l \in
\mathcal{J},
\end{equation}
and the remaining expectation $\sum_{j \in\mathcal{J}} \E
_p^{\otimes
n} [ \vert Z_j \vert^{2r} \mathbh{1} \{ \hat{j}=j \} ]$ can be
bounded by Lemma~\ref{tailsZ2} as above. Finally, the probabilities
compensate \eqref{haesslich} by Lemma~\ref{sigmabreve_to_sigmacheck}.
As concerns the expectation in \eqref{stress}, we proceed analogously using
\[
\sigma_{t,\mathrm{trunc}}^2(j \wedge\bar{j}) \leq c_{12}(
\beta,L) \bigl( \sigma_{t,\mathrm{trunc}}^2(\bar{j}) \vee
\sigma_{t,\mathrm
{trunc}}^2(j) \bigr)
\]
by Lemma~\ref{monotoneconvolution} and $\sigma_{t,\mathrm{trunc}}^2(j)
\leq c(\beta,L) \sigma_{t,\mathrm{trunc}}^2(\bar{j})$ on $A_{j,\bar
{j}} \cap\{ \hat{\sigma}_t^2(j) \leq\hat{\sigma}_t^2(\bar{j})\}$.
Combining the results for $S_1$, $S_2$ and $S_3$ proves that $R^+$ as
defined in \eqref{r+r-} is estimated by
\[
R^+ \leq c(\beta,L) \bigl( \sigma_{t,\mathrm{trunc}}^2(\bar{j}) \log n
\bigr)^{r/2}.
\]

To deduce a similar inequality for $R^-$, it remains to
investigate the probability
\[
\mathbb{P}^{\otimes n} \bigl( \hat{\sigma}_t^2(
\hat{j}) > \hat{\sigma}_t^2(\bar{j}) \bigr),
\]
since $\hat{p}_n$ and $p$ are both upper bounded by $c_1$. If $\hat
{\sigma}_t^2(\hat{j}) > \hat{\sigma}_t^2(\bar{j})$, then $\bar{j}$
cannot be an admissible exponent [see \eqref{admissiblebandwidths}],
because $\hat{j}$ had not been chosen in the minimization problem
\eqref{lepski_minimization} otherwise. Hence, by definition there
exists a multi-index $m \in\mathcal{J}$ with $\hat{\sigma}_t^2(m)
\geq
\hat{\sigma}_t^2(\bar{j})$ such that
\[
\bigl\vert\hat{p}_{n,\bar{j} \wedge m}(t) - \hat{p}_{n,m}(t) \bigr\vert>
c_3 \sqrt{\hat{\sigma}_t^2(m) \log
n}.
\]
Subsuming, we get
\begin{eqnarray*}
&&\mathbb{P}^{\otimes n} \bigl( \hat{\sigma}_t^2(
\hat{j}) > \hat{\sigma}_t^2(\bar{j}) \bigr)
\\
&&\qquad\leq\sum_{m \in\mathcal{J}} \mathbb{P}^{\otimes n} \Bigl(
\bigl\vert\hat{p}_{n,\bar{j} \wedge m}(t) - \hat{p}_{n,m}(t) \bigr\vert>
c_3 \sqrt{\hat{\sigma}_t^2(m) \log
n} , \hat{\sigma}_t^2(m) \geq\hat{
\sigma}_t^2(\bar{j}) \Bigr),
\end{eqnarray*}
and we divide the absolute value of the difference of the kernel
density estimators as in \eqref{S1_inequ} into the difference of biases
$\vert b_t(\bar{j} \wedge m) - b_t(m) \vert$ and two stochastic terms
$\vert\hat{p}_{n,\bar{j} \wedge m}(t) - \E_p^{\otimes n} \hat
{p}_{n,\bar{j} \wedge m}(t) \vert$ and $\vert\hat
{p}_{n,m}(t) -
\E_p^{\otimes n} \hat{p}_{n,m}(t) \vert$. As before,
\[
\bigl\vert b_t(\bar{j} \wedge m) - b_t(m) \bigr\vert\leq2
B_t(\bar{j}) \leq2 c_{13}(\beta,L) \sqrt{
\sigma_{t,\mathrm{trunc}}^2(\bar{j}) \log n}
\]
by Lemmas \ref{biasdiff} and \ref{lemma_hopt}, leading to the inequality
\begin{eqnarray*}
&&\mathbb{P}^{\otimes n} \bigl( \hat{\sigma}_t^2(
\hat{j}) > \hat{\sigma}_t^2(\bar{j}) \bigr)
\\
&&\qquad \leq\sum_{m \in\mathcal{J}} \mathbb{P}^{\otimes n} \Bigl(
\bigl\vert\hat{p}_{n,\bar{j} \wedge m}(t) - \E_p^{\otimes n} \hat
{p}_{n,\bar{j} \wedge m}(t) \bigr\vert+ \bigl\vert\hat{p}_{n,m}(t) - \E
_p^{\otimes n} \hat{p}_{n,m}(t)\bigr \vert
\\
&&\quad\qquad > c_3 \sqrt{\hat{\sigma}_t^2(m)
\log n} - 2 c_{13}(\beta,L) \sqrt{\sigma_{t,\mathrm{trunc}}^2(
\bar{j}) \log n} , \hat{\sigma}_t^2(m) \geq\hat{
\sigma}_t^2(\bar{j}) \Bigr)
\\
&&\qquad \leq\sum_{m \in\mathcal{J}} \bigl( \mathbb
{P}^{\otimes n}(B_{1,m}) + \mathbb{P}^{\otimes n}(B_{2,m})
\bigr)
\end{eqnarray*}
with
\begin{eqnarray*}
B_{1,m} &:=& \biggl\{ \bigl\vert\hat{p}_{n,\bar{j} \wedge m}(t) - \E
_p^{\otimes n} \hat{p}_{n,\bar{j} \wedge m}(t) \bigr\vert
\\
&&{} > \frac{1}{2} \Bigl( c_3 \sqrt{\hat{
\sigma}_t^2(m) \log n} - 2 c_{13}(\beta,L)
\sqrt{\sigma_{t,\mathrm{trunc}}^2(\bar{j}) \log n} \Bigr), \hat{
\sigma}_t^2(m) \geq\hat{\sigma}_t^2(
\bar{j}) \biggr\},
\\
B_{2,m} &:=& \biggl\{\bigl \vert\hat{p}_{n,m}(t) -
\E_p^{\otimes n} \hat{p}_{n,m}(t) \bigr\vert
\\
&&{} > \frac{1}{2} \Bigl( c_3 \sqrt{\hat{
\sigma}_t^2(m) \log n} - 2 c_{13}(\beta,L)
\sqrt{\sigma_{t,\mathrm{trunc}}^2(\bar{j}) \log n} \Bigr), \hat{
\sigma}_t^2(m) \geq\hat{\sigma}_t^2(
\bar{j}) \biggr\}. %
\end{eqnarray*}
To start with the second probability, we intersect event $B_{2,m}$ with
$A_{m,\bar{j}}$ as defined in \eqref{endspurt1}. Obviously,
\[
\mathbb{P}^{\otimes n}(B_{2,m}) \leq\mathbb{P}^{\otimes n}(B_{2,m}
\cap A_{m,\bar{j}}) + \mathbb{P}^{\otimes n} \bigl(A_{m,\bar{j}}^c
\bigr).
\]
The definition of $c_3$ and Lemma~\ref{tailsZ} allow to bound the probability
%
%
\begin{eqnarray}
\label{estimate3} \mathbb{P}^{\otimes n}(B_{2,m} \cap
A_{m,\bar{j}}) &\leq&\mathbb{P}^{\otimes n} \biggl( \frac{\vert\hat
{p}_{n,m}(t) - \E
_p^{\otimes n}
\hat{p}_{n,m}(t) \vert}{\sqrt{\sigma_{t,\mathrm{trunc}}^2(m) \log
n}} >
c_{16}(\beta,L) \biggr)
\nonumber
\\[-8pt]
\\[-8pt]
\nonumber
&\leq&2 \exp\biggl( - \frac{c_{16}(\beta,L)^2 \wedge c_{16}(\beta
,L)}{4} \log n \biggr)
\end{eqnarray}
with
%
%
\begin{equation}
\label{c16} c_{16}(\beta,L) := \biggl( \frac{c_3}{2} -
c_{13}(\beta,L) \sqrt{2 \frac{1 \vee c_1 \Vert K \Vert_2^2}{c_1 \Vert K
\Vert_2^2}} \biggr) \cdot
\sqrt{\frac{1}{2} \frac{c_1 \Vert K \Vert_2^2}{1 \vee c_1
\Vert K \Vert_2^2}}.
\end{equation}
At this point, we specify a lower bound on $c_3$. Precisely, $c_3$ has
to be chosen large enough to guarantee that
%
%
\begin{equation}
\label{c4_2} \frac{c_{16}(\beta,L)^2 \wedge c_{16}(\beta,L)}{4} \geq
\frac{r
\barbeta}{\barbeta+1} +1
\end{equation}
for any $\beta$ in the range of adaptation.
Finally, by means of Lemma~\ref{sigmabreve_to_sigmacheck},
%
%
\begin{eqnarray}
\label{estimate4}&& \mathbb{P}^{\otimes n} \bigl(A_{m,\bar{j}}^c
\bigr)\nonumber\\
&&\qquad\leq\mathbb{P}^{\otimes n} \biggl( \biggl\llvert\frac{\tilde
{\sigma}_{t,\mathrm{trunc}}^2(\bar
{j})}{\sigma_{t,\mathrm{trunc}}^2(\bar{j})} - 1
\biggr\rrvert\geq\frac
{1}{2} \biggr) + \mathbb{P}^{\otimes n} \biggl(
\biggl\llvert\frac
{\tilde
{\sigma}_{t,\mathrm{trunc}}^2(m)}{\sigma_{t,\mathrm{trunc}}^2(m)} - 1
\biggr\rrvert\geq\frac{1}{2} \biggr)
\\
&&\qquad\leq4 \exp\biggl( -\frac{3}{32 \Vert K \Vert_\mathrm{sup}^2} \log^2
n \biggr),\nonumber
\end{eqnarray}
which is of smaller order than the bound in \eqref{estimate3}.
Altogether, with this restriction on $c_3$,
\[
\mathbb{P}^{\otimes n}(B_{2,m}) \leq c(\beta,L) \bigl( \sigma
_{t,\mathrm
{trunc}}^2(\bar{j}) \log n \bigr)^{r/2}.
\]
By Lemma~\ref{monotoneconvolution}, the probability $\mathbb
{P}^{\otimes n}(B_{1,m})$ can be bounded in the same way using additionally
\[
\sigma_{t,\mathrm{trunc}}^2(\bar{j} \wedge m) \leq c_{12}(
\beta,L) \bigl( \sigma_{t,\mathrm{trunc}}^2(\bar{j}) \vee
\sigma_{t,\mathrm
{trunc}}^2(m) \bigr) = c(\beta,L) \sigma_{t,\mathrm{trunc}}^2(m),
\]
because $\sigma_{t,\mathrm{trunc}}^2(\bar{j}) \leq c(\beta,L)
\sigma
_{t,\mathrm{trunc}}^2(m)$ on the event $A_{m,\bar{j}} \cap\{ \hat
{\sigma}_t^2(m) \geq\hat{\sigma}_t^2(\bar{j}) \}$. Summarizing,
%
%
\begin{equation}
\label{P(A)} \mathbb{P}^{\otimes n} \bigl( \hat{\sigma}_t^2(
\hat{j}) > \hat{\sigma}_t^2(\bar{j}) \bigr) \leq c(
\beta,L) \bigl( \sigma_{t,\mathrm
{trunc}}^2(\bar{j}) \log n
\bigr)^{r/ 2}.
\end{equation}
Finally, by Lemma~\ref{lemma_hopt},
\begin{eqnarray*}
&&\bigl( \E_p^{\otimes n} \bigl\vert\hat{p}_{n,\hat{j}}(t) - p(t)
\bigr\vert^r \bigr)^{1/r} \\
&&\qquad\leq c(\beta,L) \biggl\{ \biggl(
\frac{\log n}{n} \biggr)^{{\barbeta}/{(\barbeta+1)}} \vee\biggl(
\frac{p(t) \log n}{n}
\biggr)^{{\barbeta}/{(2 \barbeta+ 1)}} \biggr\} \sqrt{\log n}.
\end{eqnarray*}
This completes the proof of Theorem~\ref{upperbounddensity}.
\end{pf*}

\begin{pf*}{Proof of Theorem~\ref{superefficiency}}
Before we construct the densities $p_n$ and $q_n$, we first specify
their amplitudes $\Delta_n$ and $\delta_n$ in $t$, respectively. Let
%
%
\begin{eqnarray}
\label{def_deltas} %
\Delta_n &:=&n^{-{\beta_1}/{(\beta_1+1)}} \cdot
\varrho(n),\nonumber
\\
\delta_n &:=& 4 c_4 \bigl(\beta_1^*,L_1^*,r
\bigr) \biggl( \frac{\Delta
_n}{n} \biggr)^{{\beta_1}/{(2\beta_1+1)}} (\log n)^{3/2}
\\
& =& 4 c_4 \bigl(\beta_1^*,L_1^*,r
\bigr) \Delta_n \cdot\varrho(n)^{-{(\beta_1+1)}/{(2\beta_1+1)}} (\log
n)^{3/2},
\nonumber
\end{eqnarray}
for
\[
\varrho(n) := n^{\frac{\beta_1-\beta_2}{(\beta_1+1)(\beta_2+1)}}
\]
converging\vspace*{1pt} to infinity. Note first that with this choice of $\varrho
(n)$ it holds that $\Delta_n = n^{-\beta_2/(\beta_2+1)}$,
and hence tends to zero as $n$ goes to infinity.
The amplitude $\delta_n$ is smaller than $\Delta_n$ for sufficiently
large $n$, and hence also tends to zero. Furthermore, it holds
\[
\delta_n = 4 c_4 \bigl(\beta_1^*,L_1^*,r
\bigr) \cdot n^{-\frac{\beta_2}{\beta
_2+1}} \cdot n^{\frac{\beta_2-\beta_1}{(2\beta_1+1)(\beta_2+1)}} \cdot
(\log
n)^{3/2} = o \bigl( n^{-\frac{\beta_2}{\beta_2+1}} \bigr).
\]
Denote by $K( \cdot;\beta_i), i=1,2$ the univariate, symmetric and
nonnegative functions to the H\"older exponent $\beta_i$,
respectively, as defined in the supplemental article [\citet
{supp}], Section A.4,\vspace*{1pt} normalized by appropriate choices of
$c_{17}(\beta_i)$ such that both functions integrate to one. Let
$\tilde
{L}_i = \tilde{L}_i(\beta_i), i=1,2$ be such that $K( \cdot; \beta
_i) \in\mathscr{P}_1(\beta_i,\tilde{L}_i)$. Note that $K( \cdot
;h,\beta_i) := h^{\beta_i} K(\cdot/h;\beta_i)$ has the same H\"older
regularity as $K$ [as opposed to $K_h( \cdot;\beta_i) := h^{-1}
K(\cdot/h;\beta_i)$, which has the same H\"older parameter $\beta_i$
but not necessarily the same $\tilde{L}_i$].

To ensure that $p_n(t) = \Delta_n$ we use the scaled version $K(\cdot
-t;g_{1,n},\beta_1)$ for some bandwidth $g_{1,n}$ defined below,
preserving the H\"older regularity. In order to re-establish
integrability to one, a second part is added alongside. The density
$q_n$ is then defined as $p_n$ with a perturbation added and subtracted
around $t$, that is,
\begin{eqnarray*}
p_n(x) &=& K(x-t;g_{1,n},\beta_1) +
K(x-t-g_{1,n} - g_{2,n}; g_{2,n},
\beta_1) \in\mathscr{P}_1(\beta_1,L_1),
\\
q_n(x) &=& p_n(x) - K(x-t;h_n,
\beta_2) + K(x-t-2h_n;h_n,
\beta_2) \in\mathscr{P}_1(\beta_2,L_2),
\end{eqnarray*}
with
\begin{eqnarray*}
g_{1,n}& :=& \biggl(\frac{\Delta_n}{K(0;\beta_1)} \biggr)^{{1}/{\beta
_1}},\\
g_{2,n} &:=& \bigl( 1- g_{1,n}^{\beta_1+1}
\bigr)^{{1}/{(\beta_1+1)}}, \\
 h_n &:=& \biggl( \frac{\Delta_n -\delta_n}{K(0;\beta
_2)}
\biggr)^{{1}/{\beta_2}}
\end{eqnarray*}
and suitable constants $L_1$ and $L_2$ independent of $n$. The
construction of the hypotheses is depicted in Figure~\ref
{fig:supereff}. Recall that the particular construction of $K( \cdot
;h,\beta)$ does not change the H\"older parameters and note that the
classes $\bigcup_{L>0}\mathscr{C}_c\cap\mathscr{P}_1(\beta,L)$,
$0<\beta
\leq2$, are nested ($\mathscr{C}_c$ denotes the set of continuous
functions from $\mathbb{R}$ to $\mathbb{R}$ of compact support). The
bandwidth $g_{1,n}$ tends to zero, and hence $g_{2,n}$ converges to
one. In particular, $g_{2,n}$ is positive for sufficiently large~$n$.
In turn, $h_n$ ensures that $q_n(t) = \delta_n$. Note furthermore that
$\Delta_n > \Delta_n - \delta_n$ and $K(0;\beta_1) < K(0;\beta_2)$
since the constant $c_{17}(\beta)$ is monotonously increasing in
$\beta
$ and $\beta_2 < \beta_1$. Thus, $h_n$ is smaller than $g_{1,n}$ and
consequently $q_n$ is nonnegative for sufficiently large $n$.

%
\begin{figure}[t]

\includegraphics{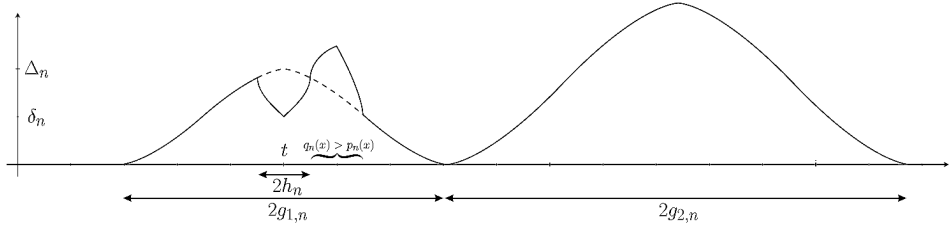}

\caption{Construction of $p_n$ (dashed line) and $q_n$ (solid line).}
\label{fig:supereff}
\end{figure}

Let $T_n(t)$ be an arbitrary estimator with property \eqref{bigrate}.
Note first that we can pass on to the consideration of the estimator
\[
\tilde{T}_n(t) := T_n(t) \cdot\mathbh{1} \bigl\{
T_n(t) \leq2 \Delta_n \bigr\},
\]
since it both improves the quality of estimation of $p_n(t)$ and
$q_n(t)$: Obviously,
\begin{eqnarray*}
\E_{p_n}^{\otimes n} \bigl\vert\tilde{T}_n(t) -
p_n(t) \bigr\vert&=& \E_{p_n}^{\otimes n} \bigl[
p_n(t) \cdot\mathbh{1} \bigl\{ T_n(t) -
p_n(t) > p_n(t) \bigr\} \bigr]
\\
&&{} + \E_{p_n}^{\otimes n} \bigl[ \bigl\vert T_n(t) -
p_n(t) \bigr\vert\cdot\mathbh{1} \bigl\{ T_n(t) -
p_n(t) \leq p_n(t) \bigr\} \bigr]
\\
&\leq&\E_{p_n}^{\otimes n} \bigl\vert T_n(t) -
p_n(t)\bigr \vert
\end{eqnarray*}
and because of $q_n(t) \leq p_n(t)$ also
\[
\E_{q_n}^{\otimes n} \bigl\vert\tilde{T}_n(t) -
q_n(t) \bigr\vert\leq\E_{q_n}^{\otimes n} \bigl\vert
T_n(t) - q_n(t) \bigr\vert.
\]
As in the proof of the constrained risk inequality in \citet
{cailowzhao2007}, by reverse triangle inequality holds
\[
\E_{q_n}^{\otimes n} \bigl\vert\tilde{T}_n(t) -
q_n(t) \bigr\vert\geq(\Delta_n - \delta_n) -
\E_{q_n}^{\otimes n}\bigl \vert\tilde{T}_n(t) -
p_n(t)\bigr \vert.
\]
In contrast to their proof, we need the decomposition:
%
%
\begin{eqnarray}
\label{supereff_decomp}&& \E_{q_n}^{\otimes n} \bigl\vert \tilde{T}_n(t)
- q_n(t)\bigr \vert
\nonumber
\\
&&\qquad\geq(\Delta_n - \delta_n) - \E_{q_n}^{\otimes n}
\bigl[ \bigl\vert T_n(t) - p_n(t) \bigr\vert\mathbh{1}_{B_n}
\bigr]
\nonumber
\\[-8pt]
\\[-8pt]
\nonumber
&&\qquad\quad{}- \E_{q_n}^{\otimes n} \bigl[ \bigl\vert\tilde{T}_n(t)
- p_n(t) \bigr\vert\mathbh{1}_{B_n^c} \bigr]
\\
&&\qquad=: (\Delta_n - \delta_n) - S_1 -
S_2,\nonumber
\end{eqnarray}
where
\[
B_n := \Biggl\{ x= (x_1, \ldots, x_n) \in
\mathbb{R}^n : \prod_{i=1}^n
\frac{q_n(x_i)}{p_n(x_i)} \leq\frac{\Delta_n}{\delta_n} \Biggr\}.
\]
By definition of $\Delta_n$ and $\delta_n$ in \eqref{def_deltas} and
the risk bound \eqref{bigrate} the first two summands in \eqref
{supereff_decomp} can be further estimated by
\begin{eqnarray*}
&&(\Delta_n - \delta_n) - S_1 \\
&&\qquad\geq(
\Delta_n - \delta_n) - \E_{p_n}^{\otimes n}
\bigl\vert T_n(t) - p_n(t) \bigr\vert\cdot\frac{\Delta
_n}{\delta_n}
\\
&&\qquad\geq(\Delta_n - \delta_n) \biggl( 1 -
\frac{c_4(\beta_1^*,L_1^*,r)
( {\Delta_n}/{n} )^{{\beta_1}/{(2\beta_1+1)}}
(\log
n)^{3/2} ({\Delta_n}/{\delta_n})}{\Delta_n - \delta_n} \biggr)
\\
&&\qquad= \delta_n \biggl( \frac{\varrho(n)^{{(\beta_1+1)}/{(2\beta_1+1)}}
(\log n)^{-3/2}}{4c_4(\beta_1^*,L_1^*,r)} - 1 \biggr)
\\
&&\qquad\quad{}\times\biggl( 1 - \frac{c_4(\beta_1^*,L_1^*,r) (
{\Delta_n}/{n} )^{{\beta_1}/{(2\beta_1+1)}} (\log n)^{3/2}
({\Delta_n}/{\delta_n})}{\Delta_n ( 1- 4c_4(\beta_1^*,L_1^*,r)
\cdot
\varrho(n)^{-{(\beta_1+1)}/{(2\beta_1+1)}} (\log n)^{3/2} )} \biggr),
\end{eqnarray*}
which is lower bounded by
\begin{eqnarray*}
(\Delta_n - \delta_n) - S_1 &\geq&
\delta_n \frac{\varrho(n)^{
{(\beta_1+1)}/{(2\beta_1+1)}} (\log n)^{-3/2}}{8c_4(\beta_1^*,L_1^*,r)}
\\
&&{}\times\biggl( 1 - \frac{2c_4(\beta_1^*,L_1^*,r) ( {\Delta_n}/{n}
)^{{\beta_1}/{(2\beta_1+1)}} (\log n)^{3/2}}{\delta_n} \biggr)
\\
&=& \delta_n \frac{\varrho(n)^{{(\beta_1+1)}/{(2\beta_1+1)}} (\log
n)^{-3/2}}{16c_4(\beta_1^*,L_1^*,r)}
\end{eqnarray*}
for sufficiently large $n$. Furthermore,
\[
S_2 \leq2 \Delta_n \cdot\mathbb{Q}_n^{\otimes n}
\bigl(B_n^c \bigr) = \delta_n
\frac{\varrho(n)^{{(\beta_1+1)}/{(2\beta_1+1)}} (\log
n)^{-3/2}}{2c_4(\beta_1^*,L_1^*,r)} \cdot\mathbb{Q}_n^{\otimes n}
\bigl(B_n^c \bigr),
\]
and it remains to show that $\mathbb{Q}_n^{\otimes n}(B_n^c)$ tends to
zero. By Markov's inequality,
\begin{eqnarray*}
\mathbb{Q}_n^{\otimes n} \bigl(B_n^c
\bigr) &=& \mathbb{Q}_n^{\otimes n} \Biggl( \prod
_{i=1}^n \frac{q_n(X_i)}{p_n(X_i)} > \frac{\Delta_n}{\delta
_n}
\Biggr)
\\
&\leq&\frac{\delta_n}{\Delta_n} \biggl( \E_{q_n} \frac
{q_n(X_1)}{p_n(X_1)}
\biggr)^n
\\
&\leq&\frac{\delta_n}{\Delta_n} \biggl( 1 + \int\frac{q_n(x)}{p_n(x)} q_n(x)
\mathbh{1} \bigl\{ q_n(x) > p_n(x) \bigr\} \,dx
\biggr)^n
\\
&\leq&\frac{\delta_n}{\Delta_n} \biggl( 1+ \frac{(2\Delta_n -
\delta
_n)^2}{K(3h_n;g_{1,n},\beta_1)} \cdot2h_n
\biggr)^n
\\
&\leq&\frac{\delta_n}{\Delta_n} \biggl( 1+ \frac{4\Delta
_n^2}{g_{1,n}^{\beta_1} K(3h_n/g_{1,n};\beta_1)} \cdot2h_n
\biggr)^n
\\
&\leq&\frac{\delta_n}{\Delta_n} \bigl( 1+ c(\beta_1,\beta_2)
\Delta_n^{{(\beta_2+1)}/{\beta_2}} \bigr)^n
\end{eqnarray*}
for sufficiently large $n$, where the last inequality is due to
\[
h_n / g_{1,n} = c(\beta_1,
\beta_2) \Delta_n^{{(\beta_1-\beta_2)}/{(\beta_1 \beta_2)}} \longrightarrow0,
\]
that is, $K(3h_n/g_{1,n};\beta_1)$ stays uniformly bounded away from
zero. Finally,
\begin{eqnarray*}
\mathbb{Q}_n^{\otimes n} \bigl(B_n^c
\bigr) &\leq&\frac{\delta_n}{\Delta_n} \exp\bigl( n \log\bigl( 1+ c(
\beta_1, \beta_2) \Delta_n^{{(\beta_2+1)}/{\beta_2}}
\bigr) \bigr)
\\
&\leq&\frac{\delta_n}{\Delta_n} \exp\bigl( n \cdot c(\beta_1,
\beta_2) \Delta_n^{{(\beta_2+1)}/{\beta_2}} \bigr)
\end{eqnarray*}
and
\[
n \Delta_n^{{(\beta_2+1)}/{\beta_2}} = 1,
\]
such that $\mathbb{Q}_n^{\otimes n}(B_n^c) \leq c(\beta_1,\beta_2)
\cdot\delta_n/\Delta_n \longrightarrow0$.
\end{pf*}

\section*{Acknowledgements}
We are very grateful to two anonymous
referees and an Associate Editor for three constructive and detailed
reports which led to a substantial improvement of our presentation and
stimulated further interesting research.


\begin{supplement}[id=suppA]
\stitle{Supplement to ``Adaptation to lowest density regions with
application to support recovery''}
\slink[doi]{10.1214/15-AOS1366SUPP} 
\sdatatype{.pdf}
\sfilename{aos1366\_supp.pdf}
\sdescription{Supplement A is organized as follows. Section A.1
contains the proofs of Lemmas \ref{varapprox}--\ref{tailsZ2},
which are central ingredients for the proof of Theorem~\ref
{upperbounddensity}. Section A.2 is concerned with the remaining proofs
of Section~\ref{sec:densitybounds}. Section A.3 contains the proofs of
Section~\ref{sec:supportapplication}. Section A.4 introduces a specific
construction of a kernel function with prescribed H\"older regularity,
which is frequently used throughout the article.}
\end{supplement}

%

%


\printaddresses
\end{document}